\documentclass[leqno]{amsart}
\usepackage{graphicx, amsmath, amssymb, amsthm, hyperref, verbatim, multicol, mathrsfs, tikz, caption, subcaption}
\usepackage{stmaryrd}
\usepackage{mathtools}  
\usepackage{cleveref}   
\usepackage{todonotes}
\usepackage{tikz-cd}
\usepackage{physics}
\usepackage{enumitem}   

\newtheorem{thm}{Theorem}[section]
\newtheorem{lemm}[thm]{Lemma}
\newtheorem{cor}[thm]{Corollary}

\newtheorem{defi}[thm]{Definition}

\newtheorem{prop}[thm]{Proposition}
\newtheorem{ques}[thm]{Question}

\theoremstyle{definition}
\newtheorem{exm}[thm]{Example}

\theoremstyle{remark}

\newcommand{\C}{\mathbb{C}}
\newcommand\numberthis{\addtocounter{equation}{1}\tag{\theequation}}

\DeclareMathOperator{\diam}{diam}

\DeclareMathOperator{\Mod}{mod}
\DeclareMathOperator{\Fill}{Fill}
\DeclareMathOperator{\Area}{Area}
\DeclareMathOperator{\apmd}{ap\,md}
\DeclareMathOperator{\ap}{ap}
\DeclareMathOperator{\J}{\mathbf{J}}

\begin{document}

\title[The branch set of minimal disks]{The branch set of minimal disks in metric spaces}

\author{Paul Creutz}
	\address{Department of Mathematics and Computer Science, University of Cologne, Weyertal 86-90, 50931 K\"oln, Germany.}
	\email{pcreutz@math.uni-koeln.de}

\author{Matthew Romney}
	\address{Mathematics Department, Stony Brook University, Stony Brook NY, 11794, USA.}
	\email{matthew.romney@stonybrook.edu}

\thanks{Both authors were supported by Deutsche Forschungsgemeinschaft grant SPP 2026. \newline {\it 2010 Mathematics Subject Classification.} Primary 49Q05. Secondary 53A10, 53C23.
	}

\begin{abstract}
We study the structure of the branch set of solutions to Plateau's problem in metric spaces satisfying a quadratic isoperimetric inequality. In our first result, we give examples of spaces with isoperimetric constant arbitrarily close to the Euclidean isoperimetric constant $(4\pi)^{-1}$ for which solutions have large branch set. This complements recent results of Lytchak--Wenger and Stadler stating, respectively, that any space with Euclidean isoperimetric constant is a CAT($0$) space and solutions to Plateau's problem in a CAT($0$) space have only isolated branch points. We also show that any planar cell-like set can appear as the branch set of a solution to Plateau's problem. These results answer two questions posed by Lytchak and Wenger. Moreover, we investigate several related questions about energy-minimizing parametrizations of metric disks: when such a map is quasisymmetric, when its branch set is empty, and when it is unique up to a conformal diffeomorphism.  
\end{abstract} 

\maketitle

\section{Introduction} \label{sec:introduction}
\subsection{Branch sets of minimal disks}
\label{sec:mindisks}
Plateau's problem asks for the existence of a disk of minimal area whose boundary is
a prescribed Jordan curve of finite length. While the problem has a longer history, the first rigorous general solutions of Plateau's problem in Euclidean space were given independently by Douglas \cite{Dou:31} and Rad\'o \cite{Rad:30} in the 1930s. This has since been generalized to various other settings, including homogeneously regular Riemannian manifolds \cite{Mor:49}, metric spaces satisfying various notions of bounded curvature~\cite{Nik79,Jos:94,MZ:10}, and certain Finsler manifolds \cite{OvdM:14,PvdM:17}. Note that arbitrary area minimizers need not have any good regularity properties. Thus Plateau's problem is usually solved by minimizing some energy functional or integrand among mappings in a suitable Sobolev space. For a metric space $X$ and Jordan curve $\Gamma \subset X$, let $\Lambda(\Gamma,X)$ denote the family of Sobolev maps~$u \in N^{1,2}(\mathbb{D}, X)$ whose  trace~$u|_{\mathbb{S}^1}$ is a weakly monotone parametrization of~$\Gamma$. We refer to a map $u \in \Lambda(\Gamma, X)$ of minimal area and minimal energy among area-minimizers as a \emph{solution to Plateau's problem}  or a \emph{minimal disk}. 
Precise definitions of the terms used here are given in \Cref{sec:definitions}.

A unified approach to the above results was given by Lytchak--Wenger in \cite{LW:17} using the quadratic isoperimetric inequality as a basic axiom. A metric space $X$ satisfies a \emph{$C$-quadratic isoperimetric inequality} if every closed Lipschitz path 
$\gamma\colon \mathbb{S}^1 \to X$ bounds a Sobolev disk $f \in N^{1,2}(\mathbb{D}, X)$ of area at most $C\cdot\ell(\gamma)^2$. More generally, $X$ satisfies a \emph{$(C,l_0)$-local quadratic isoperimetric inequality }if this property holds for all paths $\gamma$ of length at most $l_0$. In \cite{LW:17}, it is shown that if $X$ is a proper metric space satisfying a $(C,l_0)$-quadratic isoperimetric inequality, then every Jordan curve of sufficiently small length  admits a locally H\"older continuous solution to Plateau's problem that moreover extends to a continuous function on $\overline{\mathbb{D}}$. 
The assumption of a quadratic isoperimetric inequality is quite broad. It holds for all the classes of spaces mentioned above, as well as for more exotic spaces such as higher-dimensonal Heisenberg groups equipped with the Carnot--Carath\'eodory metric.  

With the existence of solutions to Plateau's problem established, one can study their topological properties. Remarkable results due to Osserman and Gulliver--Osserman--Royden show that solutions in $\mathbb{R}^3$ are immersions and hence local embeddings~\cite{Oss:70,GOR:73}. On the other hand, complex varieties are area minimizers and hence there are many examples of solutions in~$\mathbb{R}^n$ for $n\geq 4$ that are not immersed~\cite{Fed:65}. 
Nevertheless, any solution to Plateau's problem in~$\mathbb{R}^n$ is a branched immersion. This fact has recently been generalized by Stadler to solutions of Plateau's problem in $\textnormal{CAT}(0)$ spaces~\cite{Sta:ar} in the sense that any minimal disk in a CAT(0) space is a local embedding outside a finite set of branch points.

It is natural to ask how the branch set behaves for general metric spaces satisfying a quadratic isoperimetric inequality. Recall that, by the classical isoperimetric inequality, the Euclidean plane satisfies a quadratic isoperimetric inequality with constant $C = (4\pi)^{-1}$. According to a recent deep result of Lytchak--Wenger \cite{LW:18b}, any proper geodesic metric space satisfying a quadratic isoperimetric inequality with the Euclidean constant $(4\pi)^{-1}$ is a CAT($0$) space. In particular, by the above result of Stadler, the branch set of a solution to Plateau's problem is finite. The situation is different for spaces with larger isoperimetric constant. A standard example is to construct a surface by collapsing a nonseparating continuum $E$ in the Euclidean plane to a point. As discussed in Example 11.3 of \cite{LW:18a}, taking $E$ to be a closed disk or a line segment in $\overline{\mathbb{D}}$, the quotient surface $Z = \overline{\mathbb{D}}/E$ with the length metric induced by the Euclidean distance satisfies a $(2\pi)^{-1}$-isoperimetric inequality, where the constant $(2\pi)^{-1}$ is optimal.

Motivated by these examples, Lytchak and Wenger ask in Question~11.4 of \cite{LW:18a} whether there exist solutions to Plateau's problem having large branch set in metric spaces with quadratic isoperimetric constant smaller than~$(2\pi)^{-1}$. A negative answer would have been desirable for applications concerning minimal surfaces in Finsler manifolds, since every finite-dimensional normed space satisfies a $C$-quadratic isoperimetric inequality for some $C < (2\pi)^{-1}$ \cite{Cre:20}. However, it turns out that the answer to Question~11.4 in \cite{LW:18a} is affirmative, in fact for any isoperimetric constant $C> (4\pi)^{-1}$. This is our first result. 
\begin{thm} \label{thm:main}
For all $C> (4\pi)^{-1}$, there is a proper geodesic metric space $(X,d)$ satisfying a $C$-quadratic isoperimetric inequality and a rectifiable Jordan curve~$\Gamma$ in~$X$ such that any solution of Plateau's problem for $\Gamma$ collapses a nondegenerate continuum to a point.
\end{thm}
The space $X$ in Theorem~\ref{thm:main} necessarily depends on $C$, since by \cite{Cre:ar} a proper metric space satisfying a $C$-quadratic isoperimetric inequality for all $C>(4\pi)^{-1}$ is also a CAT(0) space.

The spaces we use to prove \Cref{thm:main} are a variation on the example in Section~5.1 of \cite{RRR:19}. The basic idea is to modify the quotient construction in Example~11.3 of \cite{LW:18a} by deforming the Euclidean metric on $\overline{\mathbb{D}}$ by a conformal weight that decays exponentially to zero at the collapsed set. By taking the ball of radius $\alpha>0$ centered at the point arising from the collapsed set, we obtain a family of spaces $X_\alpha$ whose isoperimetric constant goes to $(4\pi)^{-1}$ as $\alpha \to 0$. We find it interesting both that this type of construction has the properties required by \Cref{thm:main} and that an asymptotically sharp estimate of the isoperimetric constant in our family of examples can be proved. 

We note in particular that the space $X$ in \Cref{thm:main} is homeomorphic to~$\overline{\mathbb{D}}$, with $\Gamma$ the boundary of $X$. Thus our result shows that the strong geometric properties of CAT(0), and even CAT($\kappa$), surfaces are lost as soon as the Euclidean isoperimetric constant $(4\pi)^{-1}$ is replaced with anything larger. For example, CAT($\kappa$) surfaces are locally bi-Lipschitz equivalent to a disk in the Euclidean plane and thus locally Ahlfors 2-regular (see \cite[Theorem 9.10]{Res:93}), something that is not the case for our examples.\par 

In fact, there is no loss of generality to consider in Theorem~\ref{thm:main} only spaces~$X$ homeomorphic to~$\overline{\mathbb{D}}$, with $\Gamma$ the boundary curve of $X$. This is due to the following observations made in \cite{LW:18a} and \cite{Cre:19}: if $X$ satisfies a $C$-quadratic isoperimetric inequality and $u\colon \overline{\mathbb{D}}\to X$ is a solution of Plateau's problem for a given rectifiable curve $\Gamma\subset X$, then $u$ factors as \begin{equation*}
\begin{tikzcd}
\overline{\mathbb{D}} \arrow{rr}{u} \arrow[swap]{dr}{P} && X\\
& Z \arrow[swap]{ur}{\bar u}
\end{tikzcd}
\end{equation*} where
\begin{enumerate}[label=(\roman*)]
    \item \label{i1} $Z$ is a geodesic metric space homeomorphic to $\overline{\mathbb{D}}$ that satisfies a $C$-quadratic isoperimetric inequality,
    \item $P\colon \overline{\mathbb{D}}\to Z$ is a solution to Plateau's problem for $\partial Z$ in~$Z$,
    \item $\bar{u}\colon Z\to X$ is length preserving: $\ell(P\circ \gamma) = \ell(u \circ \gamma)$ for all paths $\gamma$ in $\overline{\mathbb{D}}$.
\end{enumerate}
Moreover, if $X$ itself is homeomorphic to $\overline{\mathbb{D}}$ and $\Gamma=\partial X$, then $\bar{u}$ is an isometry. In fact, Questions~11.4 and~11.5 from~\cite{LW:18a} are stated in terms of the map~$P$ rather than $u$. However, by the above observations, it suffices to consider the case that $X$ is geodesic and homeomorphic to~$\overline{\mathbb{D}}$ with~$\Gamma=\partial X$ and~$u$ is a solution to Plateau's problem for~$\Gamma$ in~$X$. For this reason, we restrict ourselves to this setting for the remainder of the paper.

\subsection{Canonical parametrizations of disks}
Throughout this section, let $Z$ be a geodesic metric space homeomorphic to~$\overline{\mathbb{D}}$ satisfying a $C$-quadratic isoperimetric inequality with rectifiable boundary $\Gamma=\partial Z$. In this situation, a map $u\in \Lambda(\Gamma,Z)$ is a minimal disk if and only if it is energy-minimizing. Thus in the following we call such a $u$ an \emph{energy-minimizing parametrization of~$Z$}.\par 

According to the classical uniformization theorem, a simply connected Riemann surface is conformally equivalent to either the Euclidean plane, the 2-sphere, or the open unit disk. The \emph{uniformization problem} asks for conditions on a metric space $Z$ homeomorphic to a given model space $X$, such as $\mathbb{R}^2$ or $\mathbb{S}^2$, that guarantee the existence of a homeomorphism from $X$ onto that space that is also geometry-preserving in some sense,~such as being quasisymmetric or quasiconformal. The prototypical result of this type is a theorem of Bonk--Kleiner that an Ahlfors 2-regular metric sphere is quasisymmetrically equivalent to the standard 2-sphere if and only if it is linearly locally connected~\cite{BK:02}. This result has since been extended in various directions, including in the papers \cite{Wil:08, Wil:10, Raj:17,Iko:19}.
    
An alternative approach to uniformization results such as the Bonk--Kleiner theorem is given in \cite{LW:20} using the existence results for Plateau's problem in metric spaces satisfying a quadratic isoperimetric inequality. In the following discussion, we limit ourselves to the case that $X = \overline{\mathbb{D}}$, although the methods can apply more generally. In \cite{LW:20}, Lytchak--Wenger show that, if $Z$ is Ahlfors 2-regular and linearly locally connected, then any energy-minimizing parametrization of $Z$ is a quasisymmetric homeomorphism. 
Moreover, such a parametrization is canonical in the sense that it is unique up to a conformal diffeomorphism of~$\overline{\mathbb{D}}$. 

In the latter part of this paper, we present several results that give a nearly complete picture of the behavior of energy-minimizing parametrizations in the quadratic isoperimetric inequality setting. First, we show that the converse implication to the theorem of Lytchak and Wenger also holds, namely that $Z$ is quasisymmetrically equivalent to $\overline{\mathbb{D}}$ only if $Z$ is Ahlfors 2-regular. Recall the standing assumption in this section that $Z$ is a geodesic metric space homeomorphic to $\overline{\mathbb{D}}$ satisfying a $C$-quadratic isoperimetric inequality with rectifiable boundary $\Gamma = \partial Z$.
\begin{thm} \label{thm:qs_doubling}
Suppose that $\Gamma=\partial Z$ is a chord-arc curve. Let $u \in \Lambda(\Gamma,Z)$ be an energy-minimizing parametrization of $Z$. The following are equivalent:
\begin{enumerate}[label=(\roman*)]
    \item $u$ is a quasisymmetric homeomorphism. \label{item:QS} 
\item $Z$ is quasisymmetrically equivalent to~$\overline{\mathbb{D}}$. \label{item:QS_equiv}
    \item $Z$ is doubling. \label{item:doubling}
    \item $Z$ is Ahlfors 2-regular. \label{item:2_reg}
\end{enumerate} 
\end{thm} 

Recall that a curve~$\Gamma$ is a called a \emph{chord-arc curve} if it bi-Lipschitz equivalent to~$\mathbb{S}^1$. Note that in \Cref{thm:qs_doubling} the disk $Z$ is not assumed  to be linearly locally connected. Rather, this follows from the quadratic isoperimetric inequality together with the chord-arc condition. We also observe that \Cref{thm:qs_doubling} does not hold without the assumption of a quadratic isoperimetric inequality. Examples of geodesic metric spaces of finite Hausdorff 2-measure quasisymmetrically equivalent to $\overline{\mathbb{D}}$ that are not Ahlfors 2-regular are constructed in \cite{Rom:19}. The assumption that $\partial Z$ is a chord-arc curve also cannot be dropped, as seen by considering the example of a Euclidean domain with outward-pointing cusp, equipped with the intrinsic metric, which is doubling but not quasisymmetrically equivalent to the unit disk. 

We now discuss assumptions under which the energy-minimizing parametrization is a quasiconformal homeomorphism in the sense of the conformal modulus definition. Note that any quasisymmetric mapping $u$ as in \Cref{thm:qs_doubling} is also quasiconformal, although the reverse is not true in general. A recent theorem of Rajala~\cite{Raj:17} characterizes metric spaces quasiconformally equivalent to a domain in $\mathbb{R}^2$ in terms of a condition on conformal modulus called \emph{reciprocality}. Roughly speaking, reciprocality requires that the moduli of path families associated to quadrilaterals and annuli are not too large; see~\cite[Definition 1.3]{Raj:17} for the precise definition. In the following result, we observe that if a reciprocal metric space also satisfies a quadratic isoperimetric inequality, then a quasiconformal parametrization is obtained from a solution to Plateau's problem, analogously to in the quasisymmetric case. 
\begin{thm}
\label{thm:rec}
If $Z$ is quasiconformally equivalent to $\overline{\mathbb{D}}$, then any energy-minimizing parametrization $u \in \Lambda(\Gamma,Z)$ is a quasiconformal homeomorphism. Moreover, $u$ is unique up to a conformal diffeomorphism of~$\overline{\mathbb{D}}$.
\end{thm}
In the situations of \Cref{thm:qs_doubling} and \Cref{thm:rec}, every energy-minimizing parametrization of $Z$ is a homeomorphism. This is not true in general, as seen from the examples discussed in \Cref{sec:mindisks}. However, by~\cite[Theorem 8.1]{LW:18a}, the energy-minimizing parametrization~$u$ is a uniform limit of homeomorphisms $\overline{\mathbb{D}} \to Z$ and thus its fibers are cell-like sets, i.e., non-separating continua. Question~11.5 in \cite{LW:18a} asks whether arbitrary cell-like subsets of $\overline{\mathbb{D}}$, such as non-contractible sets and sets that are not Lipschitz retracts, can appear as fibers of energy-minimizing parametrizations. The following result answers this question by showing that,  up to topological equivalence, arbitrary cell-like sets can occur as fibers even for a single disk~$Z$.

\begin{thm} \label{thm:other_continua}
Let $E$ be a cell-like subset of~$\mathbb{R}^2$, $z \in Z$ and $u$ be an energy-minimizing parametrization of~$Z$ such that $u^{-1}(z)$ is nondegenerate. Then there is another energy-minimizing parametrization~$v$ of~$Z$ such that $v^{-1}(z)$ is homeomorphic to~$E$.
\end{thm}
In particular, by Theorem~\ref{thm:main}, we can produce a fiber homeomorphic to $E$ in a space~$Z$ with isoperimetric constant arbitrarily close to the Euclidean constant. In fact, the branch set obtained by combining the example in Theorem~\ref{thm:main} with the proof of Theorem~\ref{thm:other_continua} is the singular fiber $v^{-1}(z)$. Thus, up to topological equivalence, we can also produce arbitrary cell-like sets as branch sets.\par 
Another consequence of Theorem~\ref{thm:other_continua} is that whenever, in a given space $X$ satisfying a quadratic isoperimetric inequality, one can rule out the existence of a minimal disk that is constant on an open ball in the interior of $\mathbb{D}$, then the fibers of any minimal disk in~$X$ must be totally disconnected. In particular, in this case any minimal disk~$u$ bounding a given Jordan curve~$\Gamma\subset X$ must map $\mathbb{S}^1$ homeomorphically onto $\Gamma$. For example, the result of Stadler discussed above rules out such collapsing in CAT(0) spaces, and hence we are able to partially recover Theorem~A in~\cite{M:00}.

In the settings of \Cref{thm:qs_doubling} and \Cref{thm:rec}, energy-minimizing parametrizations are uniquely determined up to conformal diffeomorphism of $\overline{\mathbb{D}}$.
This fact justifies the term ``canonical parametrization'' used in~\cite{LW:20} for linearly locally connected, Ahlfors 2-regular disks. However, our next result shows that uniqueness up to conformal diffeomorphism does not hold in general. This is another consequence of Theorem~\ref{thm:other_continua}.
\begin{cor} 
\label{thm:uniqueness}
If $Z$ admits a energy-minimizing parametrization that is not a homeomorphism, then there are infinitely many energy-minimizing parametrizations that are pairwise not equivalent up to a conformal diffeomorphism of~$\overline{\mathbb{D}}$.
\end{cor}
Again, by \Cref{thm:main}, such disks may be produced with arbitrary close to Euclidean isoperimetric constant.

In light of \Cref{thm:uniqueness}, one might conjecture that energy-minimizing parametrizations are unique up to a conformal diffeomorphism provided that $Z$ admits \emph{some} homeomorphic energy-minimizing parametrization. However, this is false. In \Cref{exm:homeomorphism_not_invariant}, we give an example of a disk~$Z$ for which some energy-minimizing parametrizations are homeomorphic and others are not. This is related to the topic of \emph{removable sets for conformal maps}, also known as \emph{negligible sets for extremal distance} or \emph{sets of absolute area zero}, in complex analysis; see  \cite{AB:50} and \cite{You:15}. 


Next, one might conjecture that energy-minimizing parametrizations are unique up to a conformal diffeomorphism provided that \emph{every} energy-minimizing parametrization is a homeomorphism. We expect, based on Theorem 16 in \cite{AB:50}, that the answer is also negative in general. On the other hand, we ask whether the converse to \Cref{thm:rec} holds.
\begin{ques}
Suppose that energy-minimizing parametrizations of $Z$ are unique up to a conformal diffeomorphism of $\overline{\mathbb{D}}$. Is then $Z$ quasiconformally equivalent to~$\overline{\mathbb{D}}$?
\end{ques}


 
 \subsection{Outline of the proofs}
 We explain the main ideas used to prove \Cref{thm:main}, since this is the most involved of our results. For our construction, we define a singular Riemannian metric on $\mathbb{R}^2$ by a radially symmetric conformal factor that vanishes on the unit disk and satisfies certain decay conditions, listed in \Cref{def:admissible}. Let $(X,d)$ denote the resulting metric space and $o$ the point in $X$ corresponding to the collapsed unit disk in~$\mathbb{R}^2$. We denote by $X_\alpha$ the metric space obtained by equipping the ball $\overline{B}(o,\alpha)$ with the subspace metric. Our main task is to show that the isoperimetric constant of the spaces $X_\alpha=\overline{B}(o,\alpha)$ 
 converge to~$(4\pi)^{-1}$ as~$\alpha\to 0$. The space $X$ used to prove~\Cref{thm:main} is then the space $X_\alpha$ for some sufficiently small $\alpha$. 
 
 To prove the asymptotic convergence of the isoperimetric constants, we consider $X_\alpha$ equipped with the rescaled metric $d/\alpha$. Let $\widehat{X}_\alpha$ denote the resulting metric space. Clearly $X_\alpha$ and $\widehat{X}_\alpha$ have the same isoperimetric constant. In light of Stadler's result on the branch set of minimal disks in CAT($0$) spaces, it is necessary that our construction have some region of positive curvature. In fact, the space $\widehat{X}_\alpha$ is a smooth Riemannian manifold of positive and non-integrable curvature outside the point $o$.  On the other hand, the curvature diverges to $+\infty$ at $o$ slowly enough that the rescaled spaces $\widehat{X}_\alpha$ asymptotically flatten out away from the origin and hence converge in the sense of ultralimits to a CAT(0) space.
 
 One potential approach to proving the convergence of the isoperimetric constants could be to explicitly determine isoperimetric sets in the spaces $\widehat{X}_{\alpha}$. However, even when there are natural candidates for isoperimetric sets in a given space, it is typically difficult to prove that a candidate set is indeed an isoperimetric set; see for example~\cite{MR:02,FM:13}.  
Indeed, this approach seems intractable for the spaces we construct, since balls around the origin are not isoperimetric sets and there are also no other natural candidates available. Instead, we consider a sequence of almost isoperimetric sets $A_\alpha \subset \widehat{X}_\alpha$ and show that the isoperimetric ratio $|A_\alpha|/\ell(\partial A_\alpha)^2$ converges to~$(4\pi)^{-1}$. It follows from our choice of conformal factor that, away from the origin, the spaces $\widehat{X}_\alpha$ are increasingly well approximated by a corresponding CAT(0) space, the unit ball in the cone over a circle of length~$\ell(\partial \widehat{X}_\alpha)$. The main remaining observation is that boundary length and area of almost isoperimetric sets~$A_\alpha$ cannot accumulate too much near~$o$; see~\Cref{lemm:area_bound}. One ingredient in this proof that might be useful in other situations is the general fact (\Cref{lemm:chord_arc}) that almost isoperimetric sets satisfy a coarse chord-arc type condition.

We give a short discussion of the other proofs. For the proof of~\Cref{thm:qs_doubling}, the main implication is to show that the doubling condition implies Ahlfors 2-regularity. We first prove the respective equivalence of the doubling condition and Ahlfors-$2$-regularity for proper metric spaces homeomorphic to~$\mathbb{R}^2$ in~\Cref{prop:plane}. To prove \Cref{thm:qs_doubling}, we embed the given disk isometrically in a metric plane, similarly to the constructions in~\cite{Cre:19,Cre:ar,Sta:ar}, and deduce the result for the disk from the planar one. To prove \Cref{thm:other_continua}, we apply the Riemann mapping theorem to the complement of the nondegenerate fiber to obtain additional parametrizations of $Z$. As in \cite{Oss:70},  conformal invariance of energy implies that these discontinuous parameter transforms create new energy-minimizing parametrizations.

\subsection{Organization}
The paper is organized as follows. In \Cref{sec:definitions}, we give an overview of the relevant definitions and notation related to metric Sobolev spaces, Plateau's problem, and quasiconformal mappings. \Cref{sec:chord_arc} contains the lemma on the chord-arc condition mentioned above. The proof of \Cref{thm:main} is given in \Cref{sec:main_theorem}. 
The proofs of \Cref{thm:qs_doubling} and \Cref{thm:rec} are in \Cref{sec:canonical_param}. \Cref{sec:uniqueness} contains the proofs of \Cref{thm:other_continua} and \Cref{thm:uniqueness}, as well as \Cref{exm:homeomorphism_not_invariant}.    


\section{Preliminaries} \label{sec:definitions}
 
Let $\norm{\cdot}$ denote the Euclidean norm on $\mathbb{R}^2$ and $\mathcal{L}^2$ denote Lebesgue measure on $\mathbb{R}^2$. We use $ds_{\norm{\cdot}}$ to denote the Euclidean length element. In Cartesian coordinates $(x_1,x_2)$, this is given by $ds_{\norm{\cdot}} = \sqrt{dx_1^2 + dx_2^2}$. 
For $E\subset \mathbb{R}^2$ and $x\in \mathbb{R}^2$ let $d_{\norm{\cdot}}(x,E) = \inf\{ \norm{x-y}: y\in E\}$ denote the Euclidean distance from $x$ to $E$.
Let $\mathbb{D}$ denote the open unit disk in $\mathbb{R}^2$. More generally, let $\mathbb{D}_k = \{x \in \mathbb{R}^2: \abs{x} < k\}$ denote the open disk of radius $k >0$. 

Throughout this section, $(X,d)$ and $(Y,d')$ will be metric spaces. The diameter of a set $A \subset X$ is $\diam(A) = \sup\{d(x,y): x,y\in A\}$.  For $x\in X$ and $r>0$, denote by $B(x,r)$ the open ball $\{y\in X:d(x,y)<r\}$, by $\overline{B}(x,r)$ the closed ball $\{y\in X:d(x,y)\leq r\}$ and by $S(x,r)$ the sphere $\{y\in X:d(x,y)=r\}$. Given $\alpha>0$, we use $\alpha X$ to denote the space $X$ equipped with the metric $\alpha d$.

	A \emph{path} is a continuous map $\gamma\colon I \to X$, where $I$ is a real interval or $\mathbb{S}^1$. We call $\gamma$ a \emph{closed path} in the case that $I=\mathbb{S}^1$. The length of a path $\gamma$ is denoted by $\ell(\gamma)$. The space $X$ is \emph{geodesic} if for every $x,y\in X$ there is a path $\gamma$ connecting $x$ to $y$ of length $\ell(\gamma)=d(x,y)$. Recall that, if $X$ is geodesic, then $\overline{B}(x,r)$ is the topological closure of $B(x,r)$. A \emph{curve} is the image of a path. The image of the path $\gamma$ is denoted by $\Im(\gamma)$. A set $\Gamma \subset X$ is a \emph{Jordan curve} if it is homeomorphic to $\mathbb{S}^1$. 
	

\subsection{Geometric notions}
Let $p\geq 0$.
For a metric space $X$, the \emph{Hausdorff $p$-measure} of a set $E \subset X$ is defined as
	\[\mathcal{H}^p(E) = \lim_{\varepsilon \to 0} \inf \left\{ \sum_{j=1}^\infty \omega(p) \diam(A_j)^2: E \subset \bigcup_{j=1}^\infty A_j, \diam A_j < \varepsilon \right\},  \]
	where $\omega(p) = \pi^{p/2}/\Gamma(d/2+1)$. 
	Note that the normalization constant $\omega(n)$ guarantees that the Hausdorff $n$-measure with respect to the Euclidean metric on $\mathbb{R}^n$ coincides with $n$-dimensional Lebesgue measure. In this paper, the term \emph{area} refers to the Hausdorff 2-measure and we also write $|A|$ for $\mathcal{H}^2(A)$.
	
	A metric space $X$ is \emph{Ahlfors $2$-regular} if there exists a constant $C \geq 1$ such that $C^{-1}r^2 \leq \mathcal{H}^2(B(x,r)) \leq Cr^2$ for all $x \in X$ and $r \in (0, \diam X)$. We say that $X$ is \emph{Ahlfors $2$-regular up to some scale} if there is $r_0>0$ such that this estimate holds for all $r<r_0$. A metric space $X$ is \emph{doubling} if there exists $C \geq 1$ such that for all $x \in X$ and $r>0$, the ball $B(x,2r)$ can be covered by $C$ balls of radius $r$. We say that $X$ is \emph{doubling up to some scale} if there is $r_0>0$ such that this holds for all $r<r_0$. 
	
	For a given $\lambda \geq 1$ and $R>0$, a metric space $X$ is \emph{$(\lambda, R)$-linearly locally connected} if the following two properties hold:
	\begin{enumerate}
	    \item[(LLC1)] If $x \in X$ and $r \in (0,R]$, then for all $y,z \in B(x,r)$ there exists a continuum $E \subset B(x,\lambda r)$ such that $y,z \in E$ \label{item:llc1} 
	    \item[(LLC2)] If $x \in X$ and $r \in (0,R]$, then for all $y,z \in X \setminus B(x,r)$ there exists a continuum $E\subset X\setminus B(x,r/\lambda)$ such that $y,z \in E$. \label{item:llc2}
	\end{enumerate}
	In this case, we say that $X$ is \emph{linearly locally connected up to some scale}. The space $X$ is \emph{linearly locally connected} if the above holds for some $\lambda \geq 1$ and all $R>0$. Recall that a continuum is a compact and connected set. A continuum is said to be \emph{nondegenerate} if it consists of more than one point. Note that, if $X$ is geodesic, then balls are connected and hence (LLC1) is automatically satisfied for $\lambda =1$ and any $R>0$. Also, note that a bounded connected metric space that is linearly locally connected up to some scale is in fact linearly locally connected.
	
	For a given $\lambda \geq 1$ and $R>0$, a metric space $(X,d)$ is \emph{$(\lambda, R)$-linearly locally contractible} if every ball $B(x,r) \subset X$ of radius $r\leq R$ can be contracted to a point inside the ball $B(x,\lambda r)$. In this case, we say that $X$ is \emph{linearly locally contractible up to some scale}. The space $X$ is \emph{linearly locally contractible} if this holds for all $R \leq \diam(X)/\lambda$. 
	
	The two properties of linear local connectedness and linear local contractibility are quantitatively equivalent for compact connected topological 2-manifolds \cite[Lemma 2.5]{BK:02}. Similarly, we have the following fact for metric surfaces homeomorphic to $\mathbb{R}^2$.
	
	\begin{lemm} \label{lemm:llc}
	Let $X$ be a proper metric space homeomorphic to~$\mathbb{R}^2$. 
	Then the following implications hold:
	\begin{enumerate}[label=(\roman*)]
	    \item If $X$ is $(\lambda,R)$-linearly locally contractible, then $X$ is $(\lambda',R)$-linearly locally connected for all $\lambda'>\lambda$. \label{item:contractible_llc}
	    \item If $X$ is $(\lambda,R)$-linearly locally connected, then $X$ is $(\lambda,R')$-linearly locally contractible for $R'=R/\lambda$.  \label{item:llc_contractible}
	\end{enumerate}
	\end{lemm}
	Recall that a metric space $X$ is \emph{proper} if all closed bounded subsets of $X$ are compact.
	\begin{proof}
	This follows from modifying the proof of Lemma~$2.5$ in~\cite{BK:02}, where the corresponding result is proved for closed surfaces. To obtain \ref{item:contractible_llc}, note first that the proof of (LLC1) follows exactly as in~\cite{BK:02}. To prove (LLC2), we apply properness to obtain compactness of the closed balls $K_1,K_2$ used in the proof of Lemma~$2.5$ in~\cite{BK:02}. Note here that the isomorphism $H_{1}(Z,Z\setminus K)\simeq \check{H}^{n-1}(K)$ follows as in~\cite{BK:02} from~\cite[Theorem 17, p. 296]{Spa:81} after embedding $Z$ into its one-point compactification~$\mathbb{S}^2$ and applying~\cite[Theorem 1, p. 188]{Spa:81}. The proof of \ref{item:llc_contractible} follows the corresponding proof in \cite{BK:02} without modification, even without assuming properness. Moreover, since we assume that $X$ is homeomorphic to $\mathbb{R}^2$ rather than an arbitrary closed 2-manifold as in~\cite{BK:02}, we obtain the quantitative statement given in \ref{item:llc_contractible}.
	\end{proof}
	

	\subsection{Modulus and metric Sobolev spaces}
	
	 There are various notions of Sobolev spaces with metric space targets. We use the so-called \emph{Newtonian Sobolev space}. See \cite{HKST:15} for further background and equivalence to other definitions. Let $\Delta$ be a family of paths in $X$. A Borel function $\rho\colon X \to [0, \infty]$ is \emph{admissible} for~$\Delta$ if $\int_\gamma \rho\,ds \geq 1$ for all locally rectifiable paths $\gamma \in \Delta$.  The \emph{(conformal) modulus} of $\Delta$ is defined as
	\[ \Mod \Delta = \inf \int_{X} \rho^2\,d\mathcal{H}^2, \]
	the infimum taken over all admissible functions $\rho$ for $\Delta$. A property is said to hold for \emph{almost every path} in $\Delta$ if there is a subfamily $\Delta_0 \subset \Delta$ of zero modulus such that the property holds for all $\gamma \in \Delta \setminus \Delta_0$. 
	
	Let $\Omega\subset \mathbb{R}^2$ be a bounded domain and $u\colon\Omega\to X$ be a Borel function. The function $g \colon \Omega \to [0,\infty]$ is an \emph{upper gradient} for $u$ if
	\[d((u\circ\gamma)(a), (u\circ\gamma)(b)) \leq \int_\gamma g\,ds \]
	for every rectifiable path $\gamma\colon [a,b] \to X$. We say that $u$ is a \emph{weak upper gradient} if the same holds for almost every path $\gamma$ in $X$. The \emph{Dirichlet space} $D^{1,2}(\Omega,X)$ is defined as the collection of those $u\colon \Omega \to X$ having an upper gradient in $L^2(\Omega)$. It is a standard fact that every mapping $u \in D^{1,2}(\Omega,X)$ has a \emph{minimal} weak upper gradient, denoted by $g_u$, meaning that $g_u \leq g$ almost everywhere for any weak upper gradient $g$ of $u$. 
	The \emph{Reshetnyak energy} $E_+^2(u)$ of a mapping $u \in D^{1,2}(\Omega, X)$ is defined as 
	\[E_+^2(u) = \int_{\mathbb{D}} g_u^2\, d\mathcal{L}^2.\] 
If $\varphi\colon\Omega' \to \Omega$ is a conformal diffeomorphism and $u\in D^{1,2}(\Omega,X)$ then one has $u\circ \varphi\in D^{1,2}(\Omega',X)$ and \begin{equation}
\label{eq:confinvariance}
E^2_+(u)=E^2_+(u\circ \varphi).
\end{equation} Denote by $L^2(\Omega,X)$ the set of those measurable and essentially separably valued functions $u\colon \Omega\to X$ such that the function $u^x(z)=d(x,u(z))$ lies in $L^2(\Omega)$ for some and thus any $x\in X$. The \emph{Newtonian-Sobolev space} $N^{1,2}(\Omega,X)$ is defined as the intersection $L^{2}(\Omega,X)\cap D^{1,2}(\Omega,X)$.

The \emph{approximate metric derivative at $z \in \Omega$} of a map $u\colon \Omega \to X$ is the unique seminorm $\apmd u_z$ satisfying
	\[\ap \lim_{y \to z}\frac{d(u(y), u(z)) - \apmd u_z(y-z)}{\norm{y-z}} = 0,\]
	provided this exists. Here, $\ap \lim$ denotes the \emph{approximate limit}, meaning that $z$ is a density point of a measurable set $K \subset \Omega$ for which this limit exists when restricted to $K$. If $u \in N^{1,2}(\Omega,X)$, then the approximate metric derivative is defined for almost all $z \in \Omega$; see \cite[Proposition 4.3]{LW:17}. The \emph{parametrized (Hausdorff) area} of a mapping $u \in N^{1,2}(\Omega,X)$ is
	\[\Area(u) = \int_{\Omega} \J(\apmd u_z)\, d\mathcal{L}^2(z).\]
	Here, \[\mathbf{J}(s) = \frac{\pi}{\mathcal{L}^2(B_{s})},\] where $B_{s}$ is the unit ball of the seminorm $s\colon \mathbb{R}^2 \to [0, \infty)$. Also one has 
	\[g_u(z)=\max_{v\in S^1} \,\{\apmd u_z(v)\}\]
	for almost every $z\in \Omega$. Thus \begin{equation}
	\label{eq:area-energy}
	    \Area(u)\leq E^2_+(u),
	\end{equation} 
	where equality holds precisely if $u$ is weakly conformal. Here, we say a map $u$ is \emph{weakly conformal} if for almost every $z\in \Omega$ one has $\apmd u_z=\lambda(z) \norm{\cdot}$ for some $\lambda(z)\geq 0$. We say that a map $u\colon\Omega\to X$ satisfies \emph{Lusin's condition (N)} if $\mathcal{H}^2(u(N))=0$ for every subset $N\subset \Omega$ of Lebesgue measure zero. If $u\in N^{1,2}(\Omega,X)$ satisfies Lusin's condition (N), then the area formula
	\begin{equation}
	\label{eq:areaformula}
	\Area(u)=\int_X \sharp \{z\in \Omega:u(z)=x\} \, d \mathcal{H}^2(x)
	\end{equation}
	holds; see~\cite{Kar:07}.
	
	If $\Omega$ is a Lipschitz domain, then any $u \in N^{1,2}(\Omega,X)$ uniquely determines a \emph{trace} $u|_{\partial \Omega}\colon \partial\Omega \to X$, defined for $\mathcal{H}^1$-almost every $v \in \partial \Omega$. If $\Omega=\mathbb{D}$ this is given by $u|_{\mathbb{S}^1}(v) = \lim_{t \to 1} u(tv)$ for almost every $v\in \mathbb{S}^1$.
	\subsection{The quadratic isoperimetric inequality and Plateau's problem}
	\label{subsec:QII-PP}
	Recall from the introduction that $X$ satisfies a $C$-quadratic isoperimetric inequality if every closed Lipschitz path $\gamma \colon \mathbb{S}^1\to X$ equals the trace of some disk $u\in N^{1,2}(\mathbb{D},X)$ of area at most~$C\cdot \ell(\gamma)^2$. In this case, the infimal such constant~$C$ is denoted by~$C(X)$ and called the \emph{isoperimetric constant} of~$X$. If $X$ is proper, this infimum is in fact a minimum, and hence $X$ satisfies a $C(X)$-quadratic isoperimetric inequality~\cite{Cre:ar}. By \cite{LW:18b} and \cite{Res:68}, a proper geodesic metric space $X$ satisfies a $(4\pi)^{-1}$-quadratic isoperimetric inequality if and only if it is a CAT(0) space. If $X$ is a proper geodesic metric space homeomorphic to $\overline{\mathbb{D}}$ or $\mathbb{R}^2$, then by \cite[Theorem 1.4]{LW:20} there is the following more geometric characterization: $X$ satisfies a $C$-quadratic isoperimetric inequality if and only if
	\[
 \mathcal{H}^2(U)\leq C \cdot \ell(\partial U)^2.
	\]
	for every Jordan domain $U\subset X$. 
	In the following, let $\Gamma$ be a Jordan curve in $X$. We define $\Lambda(\Gamma,X)$ as the set of those maps $u\in N^{1,2}(\mathbb{D},X)$ whose trace $u|_{\mathbb{S}^1}$ has a representative that is a monotone map $\mathbb{S}^1\to \Gamma$. The \emph{filling area} of a Jordan curve $\Gamma \subset X$ is defined as 
	\[\Fill(\Gamma) = \inf\{ \Area(u): u \in \Lambda(\Gamma, X)\}.\]
	We call $u\in \Lambda(\Gamma,X)$ a \emph{solution to the Plateau problem for $\Gamma$} in $X$ if $\Area(u)=\Fill(\Gamma)$ and the energy $E^2_+(u)$ is minimal among all such $u$. We need the following consequence of the main results in \cite{LW:17}.
	\begin{thm}
	\label{thm:plateau}
	Let $X$ be a proper metric space satisfying a $C$-quadratic isoperimetric inequality, and let $\Gamma \subset X$ be a rectifiable Jordan curve. Then there is a solution $u$ to the Plateau problem for $\Gamma$ in $X$. Any such solution has a representative that satisfies Lusin's condition (N) and extends to a continuous function on $\overline{\mathbb{D}}$.
	\end{thm}
	By choosing this continuous representative, we may thus assume that a solution to the Plateau problem is a continuous map $\overline{\mathbb{D}}\to X$ that restricts to a monotone map $\mathbb{S}^1\to \Gamma$ on the boundary. 
	
	Let $Z$ be a geodesic metric space homeomorphic to $\overline{\mathbb{D}}$ that satisfies a $C$-quadratic isoperimetric inequality and has a rectifiable boundary curve $\Gamma=\partial Z$. In this case, $u\in \Lambda(\Gamma,Z)$ is a solution to the Plateau problem if and only if $E^2_+(u)$ is minimal among all mappings in $\Lambda(\Gamma,Z)$; see the proof of Theorem 2.7 in~\cite{Cre:19}. Thus, in this case we call $u$ an \emph{energy-minimizing parametrization} of~$Z$. By Theorem~1.2 and Proposition~2.9 in~\cite{LW:20}, any energy-minimizing parametrization is a cell-like map $u\colon \overline{\mathbb{D}}\to Z$. Here, a map $u$ is called \emph{cell-like} if it is surjective and $u^{-1}(z)$ is a cell-like set for any $z\in Z$. A compact subset $K$ of an absolute neighborhood retract (such as a topological manifold) is \emph{cell-like} if $K$ is contractible to a point in any neighborhood containing it. A subset $K \subset \overline{\mathbb{D}}$ is cell-like if and only if $K$ is a nonseparating continuum which does not contain $\mathbb{S}^1$. See \cite[Sec. III.15]{Dav:86} and Section~7 in~\cite{LW:18a}. 
	
	
	\subsection{Conformal and quasiconformal mappings} \label{sec:qc_maps}
	
Depending on the situation, we denote the Riemann sphere by $\widehat{\mathbb{C}}$ or $\widehat{\mathbb{R}}^2$. It turns out that a domain $G\subsetneq \widehat{\mathbb{C}}$ is simply connected if and only if $\widehat{\mathbb{C}}\setminus G$ is a cell-like set. If $G\subset\mathbb{C}$ is a bounded simply connected domain, then, by the Riemann mapping theorem, there is a conformal diffeomorphism $f\colon\mathbb{D}\to G$. If $G$ is a Jordan domain, that is, if $\partial G$ is a Jordan curve, then Carath\'eodory's theorem states that $f$ extends to a homeomorphism $\overline{\mathbb{D}}\to \overline{G}$. 

	A homeomorphism $f\colon X \to Y$ between metric spaces of locally finite Hausdorff 2-measure is \emph{quasiconformal} if there exists $K \geq 1$ such that 
	\[K^{-1}\Mod \Delta \leq \Mod f \Delta \leq K \Mod \Delta \]
	for all path families $\Delta$ in $X$. The smallest value $K$ for which the first inequality holds is called the \emph{inner dilatation} of $f$, while the smallest value $K$ for which the second inequality holds is called the \emph{outer dilatation} of $f$. 
	
	The next definition is closely related. A homeomorphism $f \colon X \to Y$ is \emph{quasisymmetric} if there exists a homeomorphism $\eta\colon [0, \infty) \to [0,\infty)$ such that 
    \[\frac{d'(f(x),f(y))}{d'(f(x),f(z))} \leq \eta\left(\frac{d(x,y)}{d(x,z)} \right) \]
    for all triples of distinct points $x,y,z \in X$. It is straightforward to show that the doubling property is quantitatively preserved by a quasisymmetric homeomorphism, as well as the properties of linear local connectedness and contractibility.
    
    We refer the reader to \cite{Vai:71} for the basic theory of quasiconformal mappings in the planar setting. A classical result is that a homeomorphism $f\colon \mathbb{R}^2 \to \mathbb{R}^2$ is quasiconformal if and only if it is quasisymmetric. This equivalence of definitions is valid in metric spaces with well-behaved geometry \cite{HKST:01} but not in general, even for metric surfaces with locally finite Hausdorff 2-measure \cite{Raj:17,Rom:19}.  
	
	


\section{Almost isoperimetric curves}
\label{sec:chord_arc}

For this section, let $X$ be a metric space satisfying a quadratic isoperimetric inequality and $C=C(X)$ the isoperimetric constant of $X$. We say that a rectifiable Jordan curve $\Gamma\subset X$ is an 
    \emph{$\varepsilon$-isoperimetric curve} if 
    \[
    \frac{\textnormal{Fill}(\Gamma)}{\ell(\Gamma)^2} \geq \frac{C}{1+\varepsilon}.
    \]
 It follows from the proof of Lemma~1.5 in~\cite{Cre:20} that actual isoperimetric curves are chord-arc curves. Since small local perturbations of a given chord-arc curve have a controlled effect on the isoperimetric ratio but can destroy the chord-arc condition, this is not true in general for $\varepsilon$-isoperimetric curves. To overcome this problem, we make the following definition. For $\delta \in (0,1)$ and $\lambda \geq 1$, a closed rectifiable Jordan curve $\Gamma$ in~$X$ is a \emph{$(\delta,\lambda)$-chord-arc curve} if, for every~$x,y\in \Gamma$ such that the lengths~$l_1,l_2$ of the two arcs of~$\Gamma$ between~$x$ and~$y$ satisfy $\delta l_2\leq l_1 \leq l_2$, one has
    \[
    l_1\leq \lambda d(x,y).
    \]

\begin{lemm} \label{lemm:chord_arc}
Let $\lambda> 1 + \sqrt{2}$ and $0<\delta<1$. Then there is a constant $\varepsilon >0$ depending only on $\delta,\lambda$ such that every $\varepsilon$-isoperimetric curve in any geodesic metric space satisfying a quadratic isoperimetric inequality is a $(\delta,\lambda)$-chord-arc curve.
\end{lemm}
\begin{proof}
Let $K=4\lambda^{-1}+2\lambda^{-2}<2$ and
\begin{equation*}
    \varepsilon =\frac{(\delta+1)^2}{\delta^2+K \delta+1}-1.
\end{equation*}
Let $\Gamma$ be a rectifiable Jordan curve that does not satisfy a $(\delta,\lambda)$-chord-arc condition. Then there exist $x,y \in \Gamma$ such that $\delta l_2\leq l_1 \leq l_2$ and $\lambda d(x,y)<l_1$.\par 
The quadratic isoperimetric inequality implies that
\begin{align*}
    \textnormal{Fill}(\Gamma)&\leq C ((l_1+d(x,y))^2+(l_2+d(x,y))^2)\\
    &< C \left(l_1^2+\frac{2}{\lambda}l_1^2+\frac{1}{\lambda^2}l_1^2+l_2^2+\frac{2}{\lambda}l_1 l_2+ \frac{1}{\lambda ^2} l_1^2\right)\\
    \numberthis &\label{eqx1} 
    \leq C \left( l_1^2+ K l_1l_2+l_2^2\right).
\end{align*}
However,
\begin{align}
  \label{eqx2}
    \frac{l_1^2+K l_1l_2+l_2^2}{(l_1+l_2)^2}&\leq \frac{(\delta l_2)^2+K \delta l_2^2+l_2^2}{(\delta l_2+l_2)^2}=\frac{1}{1+\varepsilon}.
\end{align}
Now \eqref{eqx1} together with \eqref{eqx2} imply that~$\Gamma$ is not $\varepsilon$-isoperimetric.
\end{proof}
\section{Surfaces with almost Euclidean isoperimetric constant} \label{sec:main_theorem}

\subsection{The construction} \label{sec:construction} 
The examples used to prove \Cref{thm:main} come from the following general construction. Choose a continuum $E\subset \overline{\mathbb{D}}_K$. 
Equip $\overline{\mathbb{D}}_K$ with the singular Riemannian metric $g_E=\lambda_E^2 \cdot g_{\norm{\cdot}}$, where
\[\lambda_E(x)= \frac{e^{-1/d_{\norm{\cdot}}(x,E)}}{d_{\norm{\cdot}}(x, E)^2}\] 
for $x\notin E$ and $\lambda_E(x)=0$ for $x\in E$.
Explicitly, this means that we define the length $\ell_E$ of an absolutely continuous path $\gamma \colon I \to \overline{\mathbb{D}}_K$ by the formula
\[
\ell_E(\gamma) =  \int_I \lambda_E(\gamma(t))\cdot\norm{\gamma'(t)} \, dt,
\]
and a pseudometric $d_E$ on $\overline{\mathbb{D}}_K$ by $d_E(x,y) = \inf \ell_E(\gamma), $
where the infimum is taken over all absolutely continuous paths $\gamma$ joining $x$ to $y$. We obtain a metric space, denoted by $X_E$, by identifying the set $E$ to a point $[E]$. The preceding definition is a variation on the example in Section 5.1 of \cite{RRR:19}, where the special case $E=\{0\}$ is considered. As  noted there, the space $X_E$ is not Ahlfors $2$-regular in this case.

We consider specifically two instances of the above construction. In the first example, we choose some $K>0$ and $E = E_0 = \{0\}$. In the second example, we choose some $K>1$ and $E = E_1 = \overline{\mathbb{D}}$. We use the second of these examples to prove Theorem \ref{thm:main}, though we include the first as well since it is simpler and illustrates the main ideas involved.

We perform now the parameter transform corresponding to the exponential map at the collapsed point $[E]\in X_E$. More precisely, we represent $g_E$ in polar coordinates and then parametrize the radial directions by arc length. With the corresponding change of variables for $E=E_0,E_1$ we may represent $g_E$ as a singular Riemannian metric on $[0,R]\times \mathbb{S}^1$ for some $R>0$. Throughout this section, we denote points in $[0,R] \times \mathbb{S}^1$ in coordinates by $(r,\theta)$. For $E_0$ we obtain
\[
g_{E_0}=dr^2 + \log^2(r) r^2\,d\theta^2
\]
and similarly for $E_1$ that
\[
g_{E_1}=dr^2 + \log^2(r)(1-\log(r))^2 r^2\,d\theta^2.
\]
Thus in both cases we have represented $X_E$ as a warped product $[0,R]\times_f \mathbb{S}^1$, where respectively $f(r)=r \log(1/r)$ and $f(r)=r \log(1/r)(1+\log(1/r))$. The following definition identifies the relevant properties of these functions. 
\begin{defi} \label{def:admissible}
A continuous function $f\colon [0,R] \to [0,\infty)$ is called an \emph{admissible density} if it satisfies the following:
\begin{enumerate}[label=(\alph*)]
\item $f$ is increasing, with $f(0) = 0$. \label{item:rg_increasing}
    \item $f(r)\geq r$ for all $r\in [0,R]$.\label{item:g_decreasing}
    \item For all $r \in (0,1]$, one has $\lim_{\alpha \to 0} f(\alpha r)/(rf(\alpha)) = 1$. Moreover, this convergence is uniform with respect to $r$ on each compact subset of $(0,1]$. \label{item:g_doubling}
\end{enumerate}
\end{defi}
Recall that the standard disk $\overline{\mathbb{D}}$ may be represented as $[0,1]\times_f \mathbb{S}^1$ where $f$ is the admissible density given by $f(r)=r$. The following lemma can be shown by elementary computations, and the proof is omitted.

\begin{lemm} \label{lemm:properties_g}
For $R=e^{-1}$, it holds that the function $f\colon [0,R] \to [0, \infty)$ defined by $f(r) = r\log(1/r)$ is admissible, and for $R = e^{-(1+\sqrt{5})/2}$ that the function $f\colon [0,R] \to [0, \infty)$ defined by $f(r) = r\log(1/r)(1+\log(1/r))$ is admissible.
\end{lemm}
Fix an admissible density $f\colon[0,R]\to[0,\infty)$ and let $X^f$ be the warped product $[0,R] \times_f \mathbb{S}^1$. The point $[\{0\}\times \mathbb{S}^1]$ in $X^f$ is denoted by $o$. The definition of admissible density provides us with the following geometric properties.
\begin{lemm}
\label{lem:geoest} Let $r\in (0,R)$.
\begin{enumerate}
    \item The nearest point projection $P\colon X^f\to \overline{B}(o,r)$ is well-defined and $1$-Lipschitz.\label{item:1-lip}
    \item It holds that \[\abs{\overline{B}(o,r)}\leq r \ell(S(o,r)).\] \label{item:secest}
    \item  We have \[2\pi r\leq \ell(S(o,r))\leq C \cdot \ell(S(o,r/2))\] for some $C$ independent of~$r$. \label{item:sphcomp}
\end{enumerate}
\end{lemm}
\begin{proof}
Claim \eqref{item:1-lip} follows immediately from Condition~\ref{item:rg_increasing} in Definition~\ref{def:admissible}. By Condition~\ref{item:rg_increasing}, \[\abs{\overline{B}(o,r)} = 2\pi \int_0^r f(s)\,ds\leq 2\pi f(r) r=r\cdot \ell(S(o,r)),\] and hence we obtain~\eqref{item:secest}. The left inequality in Claim~\eqref{item:sphcomp} follows from Condition~\ref{item:g_decreasing}. The right inequality is implied by continuity of~$f$ and Condition~\ref{item:g_doubling}. 
\end{proof}
For each $\alpha \in (0,R)$, define $X_\alpha^f$ as the ball $\overline{B}(o,\alpha)\subset X^f$ endowed with the subspace metric. By Lemma~\ref{lem:geoest}.\ref{item:1-lip}, the space~$X_\alpha^f$ may be represented as the warped product $[0,\alpha]\times_f \mathbb{S}^1$ and in particular is a geodesic metric space. Our aim is to show that the quadratic isoperimetric constant of the spaces $X^f_\alpha$ converges to $(4\pi)^{-1}$ as $\alpha\to 0$. Since the isoperimetric constant is invariant under rescaling, we may equivalently consider the rescaled spaces $\widehat{X}^f_\alpha=\alpha^{-1} \cdot X_\alpha^f$. Equivalently one may represent $\widehat{X}_\alpha^f$ as the warped product $[0,1]\times_{f_\alpha}\mathbb{S}^1$, 
where $f_\alpha\colon [0,1]\to [0,\infty), f_\alpha(r)=f(\alpha r)/\alpha$ is an admissible density. These rescaled spaces $\widehat{X}^f_\alpha$ have the advantage that they do not degenerate to a point in the limit. In fact, the spaces $\widehat{X}^f_\alpha$ converge in the sense of ultralimits to a CAT(0) space. This follows from Lemma~\ref{lemm:bl_equivalence} below.

For each $\beta \in (0,\infty)$, we denote by $\mathcal{C}_\beta$ the Euclidean cone over a circle of length~$\beta$. We refer to Section 3.6.2 of~\cite{BBI:01} for the definition and basic properties of metric cones. The cone $C_\beta$ may be represented as the warped product $[0,\infty)\times_{c_\beta} \mathbb{S}^1$, where $c_\beta(r)=\beta/(2\pi)\cdot r$. We note that $\mathcal{C}_\beta$ is a CAT(0) spaces if and only if $\beta\geq 2\pi$. It follows readily that $c_\beta$ is an admissible density for all $\beta\geq 2\pi$.
\begin{lemm} \label{lemm:bl_equivalence}
For each $\alpha>0$, let $Y_\alpha$ be the unit ball in the cone over a circle of length~$\ell(\partial \widehat{X}^f_\alpha)$. Then there is a homeomorphism $F_\alpha\colon \widehat{X}_\alpha^f \to Y_\alpha$ such that:
\begin{enumerate}
    \item $d(F_\alpha(x),o)=d(x,o)$ for all $x\in \widehat{X}^f_\alpha$.
    \item For fixed $k\in (0,1)$ the restriction of $F_\alpha$ to $\widehat{X}^f_\alpha \setminus B(o,k)$ is $(1+\delta_\alpha)$-bi-Lipschitz, where $\delta_\alpha\downarrow 0$ as $\alpha \to 0$.
\end{enumerate}
\end{lemm}
\begin{proof}
Let $\beta=\ell(\partial \widehat{X}^f_\alpha)$. The homeomorphism $f_\alpha$, as a set  function, is the identity map when $Y_\alpha$ and $\widehat{X}^f_\alpha$ are represented as warped products as above. Note that $\beta=2\pi f_\alpha(1)=2\pi f(\alpha)/\alpha$. Hence, by condition~\ref{item:g_doubling} in Definition~\ref{def:admissible},
\[
\frac{f_\alpha(r)}{c_\beta(r)}=\frac{f(\alpha r)}{f(\alpha) r}\to 1
\]
as $\alpha\to 0$. This implies the claim.
\end{proof}

\subsection{Non-optimal quadratic isoperimetric inequality} \label{sec:non_optimal_qii}
As the first step towards \Cref{thm:main}, we prove directly that $X^f$ satisfies the quadratic isoperimetric inequality for some constant $C$.
\begin{lemm} \label{lemm:non_optimal_isoperimetric}
For any admissible density $f\colon [0,R]\to \mathbb{R}$, the space $X^f=[0,R]\times_f \mathbb{S}^1$ satisfies a $C$-quadratic isoperimetric inequality for some $C>0$. 
\end{lemm}
\begin{proof}
Let $\Gamma \subset X$ be a rectifiable Jordan curve bounding a closed region $A$. As discussed in Section~\ref{subsec:QII-PP}, it suffices to find a uniform upper bound on the isoperimetric ratio~$\abs{A}/\ell(\Gamma)^2$. If $A$ is far apart from~$o$ relative to the length of $\Gamma$, this uniform bound is a consequence of Lemma~\ref{lemm:bl_equivalence}. To handle the case where $A$ is close to or contains the origin, we need a more refined argument involving the geometric estimates in Lemma~\ref{lem:geoest}.

For convenience of notation, we consider $\mathbb{S}^1$ as the quotient space~$[-\pi,\pi]/\sim$ under the identification $-\pi \sim \pi$. By Proposition~3 and Remark~4 in~\cite{MHH:11}, we may assume that $A$ is symmetrized with respect to the axis $\{(r,0): r \geq 0\}$. That is, for all $r_0>0$ the set $\{ \theta \in [-\pi,\pi]: (r_0,\theta)\in A\}$ is a symmetric interval. Moreover, by passing to the convex hull,  we may assume that $A$ is totally convex within $X^f$. We split into three cases based on whether $o$ lies in the interior, the exterior or on the boundary of~$A$.

In the first case, suppose that $o\in \partial A$. Then, by convexity and symmetry, we may represent $A$ as 
\begin{equation}
\label{eq:rep}
    A=\{(r,\theta): -\pi\leq \theta \leq \pi, 0\leq r\leq \bar{r}_\theta\}
\end{equation}
where $\bar{r}\colon [-\pi,\pi]\to[0,R]$ is upper semicontinuous, axis symmetric, nonincreasing on $[0,\pi]$ and such that $\bar{r}_{-\pi}=\bar{r}_\pi=0.$ We partition $A$ into sectorial regions in the following way. Let \[A_n = \{(r,\theta) \in A: 2^{-n-1}\bar{r}_0 \leq \bar{r}_\theta \leq 2^{-n}\bar{r}_0\}\] and \[L_n = \{(r,\theta) \in \partial A: 2^{-n-1}\bar{r}_0 \leq r \leq 2^{-n}\bar{r}_0\}.\] Let $\theta_n$ denote the maximal positive value of $\theta$ such that $\bar{r}_\theta \geq 2^{-n}\bar{r}_0$. So for all $\theta \in (\theta_n, \theta_{n+1})$, we have $2^{-n-1}\bar{r}_0 \leq \bar{r}_\theta \leq 2^{-n}\bar{r}_0$. Let $\ell_n$ denote the $d_f$-length of $L_n$ and set $\tau_n = \theta_{n+1} - \theta_n$. By \Cref{lem:geoest}\eqref{item:1-lip} and the fact that $o\in \partial A$, we obtain the lower bound
\begin{equation}
\label{eq:lest}
\ell_n\geq \max \left\{ \frac{\tau_n}{\pi} \ell(S(o,2^{-n-1} \bar{r}_0)), 2^{-n} \bar{r}_0\right\}.
\end{equation}
On the other hand, Lemma~\ref{lem:geoest}\eqref{item:secest} yields 
\[
\abs{A_n}\leq \frac{\tau_n}{\pi} \abs{B(o,2^{-n}\bar{r}_0)} \leq \frac{\tau_n}{\pi} \ell(S(o,2^{-n} \bar{r}_0)) 2^{-n}\bar{r}_0.
\]
Combining these two inequalities with Lemma~\ref{lem:geoest}\eqref{item:sphcomp} gives $|A_n|\leq C \ell_n^2$, where $C$ is as in \Cref{lem:geoest}\eqref{item:sphcomp}. This gives
\[\abs{A} = \sum_{n=1}^\infty \abs{A_n} \leq C \sum_{n=1}^\infty \ell_n^2 \leq C\left(\sum_{n=1}^\infty \ell_n\right)^2 \leq C\ell(\Gamma)^2\]
as desired.

In the second case, suppose that $o$ lies in the interior of $A$. We proceed similarly to the previous case. We may again describe $A$ as in~\eqref{eq:rep} except that now we have $\bar{r}_{-\pi}=\bar{r}_\pi>0$. Thus only finitely many of the $\tau_n$ are nonzero. Let $N$ denote the largest value of $n$ such that $\tau_n>0$. For each $n<N$, the above estimates for $\ell_n$ and $\abs{A_n}$ still apply. However, for the case $n=N$, we can no longer guarantee that $\ell_N\geq 2^{-N} \bar{r}_0$. Instead, we argue as follows. First, suppose that $\tau_N\leq \pi/2$. Then
\[
|A|= \sum_{n=1}^N \abs{A_n} \leq 2 \sum_{n=1}^{N-1} \abs{A_n} \leq  2C \sum_{n=1}^{N-1} \ell_n^2 \leq 2C\left(\sum_{n=1}^{N-1} \ell_n\right)^2 \leq 2C\cdot \ell(\Gamma)^2.
\]
Next, suppose that $\tau_N\geq \pi/2$. Then the estimates~(\ref{item:secest}) and~(\ref{item:sphcomp}) from Lemma~\ref{lem:geoest} give
\begin{align*}
\abs{A_N}&\leq  \frac{\tau_N}{\pi} |B(o,2^{-N}\bar{r}_0)|\leq \frac{1}{2\pi}\cdot \frac{\tau_N}{\pi} \cdot\ell(S(o,2^{-N}\bar{r}_0))^2\\
&\leq \frac{C^2}{\pi} \left( \frac{\tau_N}{\pi} \ell(S(o,2^{-N-1}\bar{r}_0))\right)^2\leq \frac{C^2}{\pi} \ell_N^2.
\end{align*}
Estimating as before, this implies $\abs{A}\leq C_0\cdot  \ell(\Gamma)^2$ where $C_0=\max\{2C, C^2/\pi\}$.

In the final case, assume that $o \notin A$. Let $\bar{r}_0$ be the largest value of $r$ such that $(r,0)\in A$. Suppose first that $B(o,\bar{r}_0/2) \cap A \neq \emptyset$. This implies that $\ell(\partial A) \geq \bar{r}_0$. Let $\widetilde{A}$ denote the convex hull of $A$ and $\{o\}$. Then we have that $|\widetilde{A}| \geq \abs{A}$ and $\ell(\partial\widetilde{A}) \leq \ell(\partial A) + \bar{r}_0 \leq 2\ell(\partial A)$. Applying one of the previous cases gives 
\[\frac{\abs{A}}{\ell_g(\partial A)^2} \leq \frac{4|\widetilde{A}|}{\ell_g(\partial \widetilde{A})^2} \leq 4C_0.\]
The remaining possibility is that $B(o,\bar{r}_0/2) \cap A = \emptyset$. However, in this case, since $A\subset B(o,\bar{r}_0)$, by Lemma~\ref{lemm:bl_equivalence} and Lemma~\ref{lem:geoest}(\ref{item:sphcomp}), the isoperimetric ratio of $A$ is comparable to the isoperimetric ratio of a set $\hat{A}$ in the cone $\mathcal{C}_\beta$ over a circle of length $\beta=\ell(S(o,\bar{r}_0))/\bar{r}_0\geq 2\pi$. Hence, since $\mathcal{C}_\beta$ is a CAT(0) space, 
\[
\frac{|A|}{\ell(\partial A)^2}\leq \frac{|\widetilde{A}|}{\ell(\partial \widetilde{A})^2}\leq \frac{(1+\delta_R)^4}{4\pi}
\]
where $\delta_R$ is chosen as in Lemma~\ref{lemm:bl_equivalence} for~$k=1/2$. Thus $X^f$ satisfies a $C_1$-quadratic isoperimetric inequality where $C_1=\max\{4C_0,(1+\delta_R)^4/(4\pi)\}$.
\end{proof}

\subsection{Asymptotically sharp isoperimetric inequality}
\label{sec:blow_ups}

Now we prove the convergence of the isoperimetric constants.
\begin{thm}
\label{thm:optima_constant}
Let $f\colon [0,R]\to \mathbb{R}$ be an admissible density, and for $\alpha\in (0,R)$ let~$X^f_\alpha$ be the warped product $[0,\alpha]\times_f \mathbb{S}^1$. Then for any $C>(4\pi)^{-1}$ there is $\alpha>0$ such that $X_\alpha^f$ satisfies a $C$-quadratic isoperimetric inequality.
\end{thm}
\begin{proof}
Let $C_\alpha$ be the isoperimetric constant of $X_\alpha^f$. By Lemma~\ref{lem:geoest}(\ref{item:1-lip}), the constants $C_\alpha$ are nondecreasing in~$\alpha$ and by Lemma~\ref{lemm:non_optimal_isoperimetric} they are finite. Our objective is to show that
\begin{equation}
\label{eq:asy}
    \lim_{\alpha\to 0} C_{\alpha}=\frac{1}{4\pi}.
\end{equation}
Since $C_\alpha$ is nondecreasing it suffices to show that $\limsup_{n\to \infty} C_{\alpha_n}=(4\pi)^{-1}$ for some sequence $(\alpha_n)$ with $\alpha_n \to 0$. Choose such a sequence~$(\alpha_n)$. We write $X_n$ for the rescaled space $\widehat{X}_{\alpha_n}^f=\alpha_n^{-1}\cdot X_{\alpha_n}^f$ as in \Cref{sec:construction}.

Fix a sequence $(\varepsilon_n)$ such that $\varepsilon_n\to 0$. For each $n \in \mathbb{N}$ we may choose an $\varepsilon_n$-isoperimetric domain $A_n\subset X_n$,~see \Cref{sec:chord_arc} for the definition. Replacing $\alpha_n$ by a smaller value if necessary, we assume that there is $x\in \partial A_n$ such that $d(o,x)=1$. We also assume that $A_n$ is symmetrized as in the proof of \Cref{lemm:non_optimal_isoperimetric}. 
Recall that, by Lemma~\ref{lemm:chord_arc}, $\partial A_n$ satisfies a $(\delta_n,\lambda)$-chord arc condition, where $\lambda=3$ and $\delta_n\to 0$. The chord-arc condition, together with the geometric estimates in \Cref{lem:geoest}, implies that the area of the sets $A_n$ cannot asymptotically concentrate at the vertex point $o$. The following lemma, whose proof is deferred to the end of this section, makes this precise.

\begin{lemm} \label{lemm:area_bound}
For each $\nu \in (0,1)$, define $b^\nu\colon \mathbb{N}\to \mathbb{R}$ by $b^\nu_n=\abs{A_n \cap B(o,\nu)}$. Then $\norm{b^\nu}_\infty \to 0$ as $\nu \to 0$.
\end{lemm}

The proof of \eqref{eq:asy} is a straightforward argument using Lemma~\ref{lemm:area_bound}, the comparison result Lemma~\ref{lemm:bl_equivalence} and the observation that Euclidean cones satisfy a~$(4\pi)^{-1}$-quadratic isoperimetric inequality.

Fix $\mu >0$. Let $\nu\in (0,\tfrac{1}{2})$ be such that $\max\{\nu, \norm{b^\nu}_\infty\} \leq \mu$. Let $N\in \mathbb{N}$ be large enough so that $F_n:=F_{\alpha_n}$ as in \Cref{lemm:bl_equivalence} is $(1+\mu)$-bi-Lipschitz on $X_n \setminus B(o,\nu)$ for~$n\geq N$. 
If $A_n$ does not intersect $B(o,\nu)$, denote by~$K^\nu_n$ the image set $F_n(A_n)$ in $Y_n$. Otherwise, denote by $K^\nu_n$ the convex hull of $F_n(A_n)\setminus B(o,\nu)$ and the vertex point $o$. In the former case we immediately have $\ell(\partial K^\nu_n)\leq (1+\mu)\ell(\partial A_n)$ and in the latter one $\ell(\partial A_n)\geq 1$, which implies
\[
\ell(\partial K^\nu_n)\leq (1+\mu) \ell(\partial  A_n)+2\nu\leq (1+3\mu) \ell(\partial A_n).
\]
Applying Lemma~\ref{lemm:area_bound}, we obtain the estimate
\[
\abs{A_n}\leq (1+\mu)^2 \abs{K_n^\nu}+b_n^\nu\leq (1+\mu)^2 \abs{K_n^\nu}+\mu.
\]
Combining these observations with the fact that $Y_n$ is a $\textnormal{CAT}(0)$ space, we obtain
\begin{equation*}
    \frac{\abs{A_n}}{\ell(\partial A_n)^2}  \leq (1+\mu)^2(1+3\mu)^2\frac{\abs{K^\nu_n}}{\ell(\partial K^\nu_n)^2}+\mu \leq \frac{(1+\mu)^2(1+3\mu)^2}{4\pi}+\mu.
\end{equation*}
Thus 
\begin{equation*}
    \limsup_{n\to \infty} C_{n}=\limsup_{n \to \infty}\frac{\abs{A_n}}{\ell(\partial A_n)^2}\leq \frac{(1+\mu)^2(1+3\mu)^2}{4\pi}+\mu.
\end{equation*}
Since $\mu >0$ is arbitrary, this completes the proof of \eqref{eq:asy}.\end{proof}

The proof of \Cref{lemm:area_bound} requires the following auxiliary lemma, stating that the boundary lengths of $A_n$ cannot accumulate too much at the origin.

\begin{lemm}
\label{lemm:lbound}
For each $\nu \in (0,1)$, define $l^\nu\colon \mathbb{N} \to \mathbb{R}$ by~$l_n^\nu=\ell(\partial A_n \cap B(o,\nu))$. Then 
$\norm{\ell^\nu}_\infty \to 0$ as $\nu\to 0$.
\end{lemm}
\begin{proof}
Fix $\mu >0$ and choose $\delta,\nu_0 >0$ accordingly to be determined later. Let $N \in \mathbb{N}$ be such that $\delta_n \leq \delta$ for all $n \geq N$.  As a first observation, $\diam (X_n)\leq 2$ and the chord-arc condition imply that $\ell(\partial A_n)\leq 2\lambda$.
We claim that \begin{equation}\label{eq:laim}
l^\nu_n\leq \mu.
\end{equation}
for all $n \geq N$ and $\nu \leq \nu_0$. 
If $A_n$ does not intersect $B(o,\nu)$ this is obvious. So lets assume otherwise. Then $\partial A_n$ divides into two connected arcs, one contained in $B(o,\nu)$ and one contained in $X_n \setminus B(o,\nu)$. We set $l_1 =\ell^\nu_n$ and $l_2=\ell(\partial A_n \setminus B(o,\nu))$.  Since $A_n$ intersects $B(o,\nu)$ and $S(o,1)$, we must have $l_2\geq 1-\nu$. Now we repeatedly apply the $(\delta_n,\lambda)$-chord-arc condition. In the case that $l_1\leq l_2$, we have
\[
l_1\leq 2\lambda \cdot \max\{\delta,\nu\}
\]
In the case that $l_2\leq l_1$,  either
\[
l_1\leq \frac{l_2}{\delta} \leq \frac{2 \lambda \nu}{\delta}
\]
or \[
1-\nu \leq l_2< \delta\cdot \ell(\partial A_n)\leq 2\lambda \delta.
\] However, the latter is contradiction for $\nu_0$ and $\delta$ sufficiently small. Thus \eqref{eq:laim} follows upon choosing a small $\delta \leq \mu/(2\lambda)$ and setting $\nu_0=\mu \delta/(2\lambda)$.

Thus we have $\norm{\ell^\nu}_\infty \leq \max\{ \mu, \ell_1^\nu,\ldots, \ell_{N-1}^\nu\}$ for all $\nu \leq \mu^2/(2\lambda L)$. 
Since $\ell^\nu\to 0$ pointwise, this establishes the lemma.
\end{proof}

\begin{proof}[Proof of \Cref{lemm:area_bound}]
Let $n\in \mathbb{N}$ and $\nu \in (0,1)$. If $A_n$ does not intersect $B(o,\nu)$, then $b^\nu_n=0$. So assume otherwise. Denote by $K_n^\nu$  the convex hull of $A_n \cap \partial B(o,\nu) $ and set $H_n^\nu =A_n\cap B(o,\nu)\setminus K_n^\nu$. Then $\ell(\partial H_n^\nu)\leq l^\nu+2\nu$. Moreover, observe that $S_n^\nu$ is contained in a sector in $B(p,\nu)$. Thus \Cref{lem:geoest}\eqref{item:secest} implies that $\abs{S_n^\nu} \leq \nu \ell(\partial A_n)\leq 2\lambda \nu$ and hence
\begin{equation}
\label{equ:gb}
b^\nu_n\leq \abs{S_n^\nu}+\abs{H_n^\nu}\leq \nu L + C (\norm{l^\nu}_\infty+2\nu)^2.
\end{equation}
where $C$ is as in Lemma~\ref{lemm:non_optimal_isoperimetric}. Observe that the right-hand side of~\eqref{equ:gb} is independent of $n$ and, by Lemma~\ref{lemm:lbound}, converges to~$0$ as~$\nu \to 0$.
\end{proof}

\subsection{Proof of \Cref{thm:main}} 
Now we complete the proof of Theorem~\ref{thm:main}. Let $C>(4\pi)^{-1}$. We consider the space $X_E$ as in \Cref{sec:construction}, taking $E = \overline{\mathbb{D}}$ and some $K>1$. By Theorem~\ref{thm:optima_constant}, for $K>1$ sufficiently small the space $X_E$ satisfies a $C$-quadratic isoperimetric inequality. Thus it suffices to show that any energy-minimizing parametrization~$u\colon \overline{\mathbb{D}}\to X_E$ collapses a nondegenerate continuum in $\overline{\mathbb{D}}$ to a point. Since any energy-minimizing parametrization is a cell-like map, it is enough to show that there is no homeomorphic energy-minimizing parametrization of $X_E$.

Let $P \colon \overline{\mathbb{D}}\to X_E$ be the canonical map given by $P(z)=[Kz]$. Then the map $P$ is continuous and is contained in $N^{1,2}(\mathbb{D},X_E)$, since it is the postcomposition of a diffeomorphism by a Lipschitz map. Thus $P\in \Lambda(\partial X_E,X_E)$. Also $P$ restricts to a conformal diffeomorphism on $A=\mathbb{D}\setminus \overline{\mathbb{D}}_r$, where $r=1/K$, and maps $\overline{\mathbb{D}}_r$ to the collapsed point $o$. The approximate metric derivative of~$P$ is given by
\[
\apmd P_z(v)=K \lambda_E(Kz) \norm{v}.
\]
This observation, together with the equality case in~\eqref{eq:area-energy} and the area formula~\eqref{eq:areaformula}, implies
\[
E^2_+(P)=\Area(P)=\mathcal{H}^2(X_E).
\]

Suppose that $u \colon \overline{\mathbb{D}}\to X_E$ is a homeomorphic energy-minimizing parametrization. Then, by Theorem~\ref{thm:plateau}, $u$ satisfies Lusin's condition (N). Hence, inequality~\eqref{eq:area-energy} and the area formula~\eqref{eq:areaformula} give
\[
E^2_+(P)=\mathcal{H}^2(X_E)\leq\Area(u)\leq E^2_+(u)\leq E^2_+(P).
\]
Thus by the equality case in \eqref{eq:area-energy} the map $u$ is weakly conformal. In particular $u$ is also a minimizer of the Dirichlet energy, which in the metric context is also called the Korevaar--Schoen energy; cf.~\cite{LW:17b}. Thus, for $p=u^{-1}(o)$ the restriction of $u$ to $\mathbb{D}\setminus \{p\}$ is a continuous and weakly harmonic map in the classical sense 
taking values in the smooth Riemannian manifold $X_E\setminus \{o\}$. Such continuous, weakly harmonic maps between Riemannian manifolds are however automatically smooth; see Theorem~9.4.1 in~\cite{Jos:17}. Thus, the composition $P^{-1}\circ u$ is a smooth weakly conformal homeomorphism $\mathbb{D}\setminus \{p\}\to A$ and hence either holomorphic or anti-holomorphic. However, such a map cannot exist by Riemann's theorem on the removability of singularities of holomorphic maps.

\section{Canonical parametrizations} \label{sec:canonical_param}

  \subsection{Quasisymmetric parametrizations} \label{sec:doubling_equivalence}
  
  In this section, we give the proof of \Cref{thm:qs_doubling}. This proof depends on the following proposition. 
  
  \begin{prop}
  \label{prop:plane}
  Let $X$ a be proper geodesic metric space homeomorphic to~$\mathbb{R}^2$ such that $\mathcal{H}^2_X$ is locally finite.
\begin{enumerate}
    \item \label{it2}
    Assume $X$ satisfies a quadratic isoperimetric inequality. Then $X$ is doubling if and only if $X$ is Ahlfors $2$-regular.
    \item \label{it3}
    Assume $X$ satisfies a local quadratic isoperimetric inequality. Then $X$ is doubling up to some scale if and only if $X$ is Ahlfors $2$-regular up to some scale. 
\end{enumerate}
  \end{prop}
  In the following proof, we say that a closed path $\gamma\colon \mathbb{S}^1 \to X$ winds around $p\in X$ if $p \notin \Im(\gamma)$ and $\gamma$ is noncontractible within $X\setminus \{p\}$. We use the following observation: if $\gamma$ winds around $p\in X$ and $\eta\colon [0,1] \to X\setminus\{p\}$ is a path connecting points $\gamma(t)$ and $\gamma(s)$, then the concatenation of $\eta$ with one of the two arcs of $\gamma$ between $\gamma(t)$ and $\gamma(s)$ must wind around~$p$.
  \begin{proof}
\eqref{it2} It is clear that any Ahlfors $2$-regular space is doubling. It remains to prove the converse. Observe that, by the proof of Theorem~$8.5$ in~\cite{LW:18a}, the assumptions that $X$ is proper, geodesic, homeomorphic to~$\mathbb{R}^2$, of locally finite Hausdorff 2-measure and satisfies a quadratic isoperimetric inequality imply that $X$ is lower Ahlfors $2$-regular. Note here that \cite[Lemma~$6.11$]{LW:18a} applied in the proof holds for any space of topological dimension~$2$ by~\cite[Theorem~VII.$2$]{HW:41}. \par 
Assume now that $X$ is also doubling. We claim that $X$ is upper Ahlfors $2$-regular. 
By the coarea inequality~\cite[2.10.25-26]{Fed:69} and because $\mathcal{H}^2_X$ is locally finite, $\mathcal{H}^1(S(x,s))$ is finite for almost every $s\in (0,\infty)$. Choose such $s\in (3r,4r)$. Since $B(x,r)$ is connected and contained in a bounded component of $X\setminus S(x,s)$, there is a connected component $K$ of $S(x,s)$ such that $x$ is contained in a bounded component of $X\setminus K$. Now $K$ is a Peano continuum by Corollary 2B in~\cite{Fre:92}. Hence, by~\cite[Theorem IV.6.7]{Wil:49}, there is a Jordan curve $\Gamma\subset K$ of finite length that winds around $x$.\par 
Since $X$ is doubling, there are a number $M$ which only depends on the doubling constant of $X$ and points $z_1,\ldots,z_M \in X$ such that $\overline{B}(x,4r)\subset \overline{B}(z_1,r)\cup \cdots \cup \overline{B}(z_M,r)$. Thus there are points $y_1,\ldots ,y_M\in \Gamma$ such that $\Gamma\subset Y$, where $Y = \overline{B}(y_1,2r)\cup \cdots \cup \overline{B}(y_M,2r)$. By the Arzel\`a--Ascoli theorem and lower semicontinuity of length, there exists a rectifiable path $\gamma\colon \mathbb{S}^1 \to Y$ that winds around $x$ and has minimal length among all paths with image in~$Y$ and winding around $x$. Since $\gamma$ has minimal length, $\gamma$ must be injective, and hence its image is a Jordan curve.  
 
 
 Next, let $z,w\in \gamma$ be such that both arcs $\gamma_+,\gamma_-$ of $\gamma$ connecting $z$ to $w$ have equal length. Since $X$ is geodesic there is a path $\eta$ connecting $z$ to $w$ within~$Y$ satisfying $\ell(\eta)\leq 4Mr$. 
 By the observation preceding the proof, either the concatenation $\gamma_+\cdot \eta$ or the concatenation $\gamma_- \cdot \eta$ must wind around $x$. Let $\nu$ be the respective path winding around $x$. Then the length minimality of $\gamma$ gives
 \begin{equation*}
 2\ell(\gamma_+)= \ell(\gamma)\leq \ell(\nu)=  \ell(\gamma_+)+\ell(\eta),
 \end{equation*}
 which implies that $\ell(\gamma_+)\leq \ell(\eta)$ and hence $\ell(\gamma)\leq 8 Mr$. Since $B(x,r)$ is connected, $\gamma$ winds around $x$ and, since $B(x,r)\cap \gamma=\emptyset$, the ball $B(x,r)$ is contained in the Jordan domain bounded by~$\gamma$. Thus the quadratic isoperimetric inequality implies
 \begin{equation*}
     \mathcal{H}^2(B(x,r))\leq C \cdot \ell(\gamma)^2\leq 4K\cdot \ell(\eta)^2\leq 16M^2C \cdot r^2,
 \end{equation*}
 where $C$ is the quadratic isoperimetric constant of $X$. This completes the proof.
\eqref{it3} This follows from the proofs of~\eqref{it2}, since both implications come with control on the scales.
\end{proof}
  
  
  \begin{proof}[Proof of Theorem~\ref{thm:qs_doubling}]
  
  \ref{item:QS} $\Longrightarrow$ \ref{item:QS_equiv}: This is immediate.
  
  \ref{item:QS_equiv} $\Longrightarrow$ \ref{item:doubling}:
It is a standard fact that the doubling property is preserved by quasisymmetric mappings and that the closed Euclidean disk $\overline{\mathbb{D}}$ is doubling. Thus, if $Z$ is quasisymmetrically equivalent to $\overline{\mathbb{D}}$, then $Z$ is doubling.

\ref{item:doubling} $\Longrightarrow$ \ref{item:2_reg}: By \cite[Theorem~1.6]{LW:18a} and the finiteness of $\mathcal{H}^2(Z)$, the space $Z$ is lower Ahlfors $2$-regular. To prove upper Ahlfors regularity, we embed~$Z$ into a metric plane~$X_Z$ and apply Proposition~\ref{prop:plane}. The construction of the plane $X_Z$ is essentially borrowed from~\cite{Cre:19}. 
Let $\partial Z$ be endowed with its intrinsic metric and let $W=\partial Z\times [0,\infty)$ with the product metric. Then we let $X_Z$ be the metric quotient of $Z\sqcup W$ under identification of~$\partial Z$ and $\partial Z\times \{0\}$. By \cite[Section 3]{Cre:19}, $X_Z$ is proper, geodesic, homeomorphic to~$\mathbb{R}^2$, contains $X$ isometrically and satisfies a local quadratic isoperimetric inequality. Thus by Proposition~\ref{prop:plane} it suffices to show that $X_Z$ is doubling.\par
Consider a ball $B(x,2r) \subset X_Z$. If $B(x,2r)\cap Z=\emptyset$, then $B(x,2r)$ is contained in $C_1$ balls $B(y_i,r)$, where $y_i\in W$ and $C_1$ is the doubling constant of~$W$. Similarly the desired estimate follows if $B(x,2r)\subset Z$. Thus assume $y\in \partial Z\cap B(x,2r)$. Then $B(x,2r)$ is contained in $B(y,4r)$. This in turn is, as a set, contained in the union of $B_Z(y,4r)$ and $B_W(y,(4\lambda+1)r)$; see~\cite[Section 3]{Cre:19}. Thus the doubling estimate follows since $Z$ and $W$ are doubling.

\ref{item:2_reg} $\Longrightarrow$ \ref{item:QS}:
We claim that if $Z$ is Ahlfors 2-regular, then the energy-minimizing parametrization $u$ is quasisymmetric. By Theorem 1.1 in~\cite{LW:20}, it suffices to show that $Z$ is linearly locally connected. The proof of this does not rely on the Ahlfors 2-regularity of $Z$ but only on the quadratic isoperimetric inequality and the chord-arc condition. 

Observe that $Z$ is linearly locally contractible by Theorem 8.6 in~\cite{LW:18a}. Thus Theorem 3.2 of \cite{Cre:19} implies that $X_Z$ is linearly locally contractible up to some scale, and hence \Cref{prop:plane} that $X_Z$ is linearly locally connected up to some scale. Let $\lambda \geq 1$, $R>0$ be such that $X_Z$ is $(\lambda,R)$-linearly locally connected. We claim that $Z$ is $(\lambda',R)$-linearly locally connected where $K$ is the chord-arc constant of $\partial Z$ and $\lambda'=2\lambda(K+1)$. This completes the proof since $Z$ is bounded. The constant $\lambda'$ is chosen such that when $w_1,w_2\in \partial Z \cap B(x,r/\lambda')$ where $x\in Z$ and $0<r\leq R$ then one of the arcs of $\partial Z$ between $w_1$ and $w_2$ is contained in $B(x,r/\lambda)$.

Note that the (LLC1) property follows immediately from $Z$ being a geodesic metric space. Thus it remains to show the (LLC2) property. To this end, let $y_1,y_2\in Z\setminus B(x,r)$ where $x\in Z$ and $0<r\leq R$. Then there is a continuum $K\subset X_Z\setminus B(x,r/\lambda) $ containing $y_1$ and $y_2$. For $i=1,2$ let $K_i$ be the connected of $K\cap Z$ in which $y_i$ lies. If $K_1=K_2$ then we are done since $\lambda'\geq \lambda$. Otherwise let $z_i\in K_i \cap \partial Z$. Since $z_i\in \partial Z\setminus B(x,r/\lambda)$ and by our choice of $\lambda'$, one of the arcs of $\partial Z$ between $z_1$ and $z_2$ must be contained in $Z\setminus B(x,r/\lambda')$. Joining $K_1$, $K_2$ and this arc gives the desired continuum. 
\end{proof}
\subsection{Quasiconformal parametrizations}

We prove \Cref{thm:rec}, which states that an energy-minimizing parametrization of $Z$ is a quasiconformal homeomorphism if $Z$ is reciprocal (Definition 1.3 in \cite{Raj:17}). By \cite[Thm. 1.4]{Raj:17}, this condition is equivalent to $Z$ admitting some quasiconformal parametrization by $\overline{\mathbb{D}}$.
\begin{proof}[Proof of \Cref{thm:rec}]
 Let $u\colon\mathbb{D}\to Z$ be an energy-minimizing parametrization. By assumption, there exists a $K$-quasiconformal mapping $g\colon \overline{\mathbb{D}} \to Z$. Then the composition $h = g^{-1} \circ u$ is a continuous monotone mapping of $\overline{\mathbb{D}}$ to itself.  We note that $u$ satisfies the modulus inequality 
$\Mod \Gamma \leq Q\cdot \Mod u \Gamma$
for all path families $\Gamma$ in $\overline{\mathbb{D}}$ for $Q=4/\pi$; see Section~3 in~\cite{LW:20}. Thus $h\colon \overline{\mathbb{D}} \to \overline{\mathbb{D}}$ satisfies
\begin{equation} \label{equ:modulus}
\Mod \Gamma \leq KQ \cdot \Mod h \Gamma
\end{equation} 
for all path families $\Gamma$ in $\overline{\mathbb{D}}$. It is standard to show that the modulus inequality \eqref{equ:modulus} and the monotonicity of $h$ imply that $h$ is a homeomorphism. For example, see Theorem 3.6 in \cite{LW:20} for a similar argument. Indeed, suppose there exists a nondegenerate continuum $E \subset \overline{\mathbb{D}}$ such that $h(E)$ is a point. Choose a nondegenerate continuum $F \subset \partial \overline{\mathbb{D}} \setminus E$. As shown in Proposition 3.5 of \cite{Raj:17}, the family of paths $\Gamma(E,F)$ that intersect both $E$ and $F$ has positive modulus. However, $\Mod h\Gamma(E,F)=0$, since the modulus of the family of all paths intersecting a point in $\overline{\mathbb{D}}$ is zero. This contradicts \eqref{equ:modulus}, and we conclude that $h$ is a homeomorphism.

Since $h$ is a homeomorphism satisfying the one-sided modulus inequality \eqref{equ:modulus}, it follows from the classical theory of planar quasiconformal mappings that $h$ is in fact quasiconformal; see Theorem 22.3 and Theorem 34.1 of \cite{Vai:71}. Thus $u$ is quasiconformal, being the composition of quasiconformal homeomorphisms.

It remains to show uniqueness up to conformal diffeomorphism. To this end, let $u$ and $v$ be energy-mininimizing parametrizations. Then $\varphi=v^{-1}\circ u$ defines a quasiconformal homeomorphism $\overline{\mathbb{D}}\to \overline{\mathbb{D}}$, and hence $\varphi$ is also quasisymmetric. From this point on, one can argue exactly as in the proof of Theorem~6.1 in~\cite{LW:20} to see that $\varphi$ is in fact a conformal diffeomorphism. Since $u=v\circ \varphi$ this completes the proof.
\end{proof}

In fact, the inverse of the map $f_0$ produced in Theorem 1.5 of \cite{Raj:17} is an energy-minimizing parametrization of $Z$. This is a consequence of the fact, shown in the proof of \cite[Thm. 1.5]{Raj:17}, that~$f_0^{-1}$ is infinitesimally isotropic in the sense of Definition 3.3 in \cite{LW:20}.  
In particular, by \cite[Thm. 1.5]{Raj:17} and its proof, any energy-minimizing parametrization $u$ has minimal outer dilatation among all quasiconformal mappings $g\colon \overline{\mathbb{D}}\to Z$ and has inner dilatation at most $2$.

\section{Surgering energy-minimizing parametrizations} \label{sec:uniqueness}

In this section, we prove \Cref{thm:other_continua} and~\Cref{thm:uniqueness} and present Example~\ref{exm:homeomorphism_not_invariant}. For all these, we will surger a given energy-minimizing parametrization to produce new ones with desired properties. The main tools are the Riemann mapping theorem, invariance of energy under precomposition by conformal diffeomorphisms and the following simple negligibility lemma.
\begin{lemm}
\label{lemm:surg}
Let $X$ be a complete metric space, $K\subset \overline{\mathbb{D}}$ be compact and $\Omega=\mathbb{D}\setminus K$. Suppose that $u\colon \overline{\mathbb{D}}\to X$ is continuous and satisfies $\mathcal{H}^1(u(K))=0$ and $u|_\Omega \in N^{1,2}(\Omega,X)$. Then $u\in N^{1,2}(\mathbb{D},X)$ and $E^2_+(u)=E^2_+(u|_\Omega)$.
\end{lemm}
With the exception of Example~\ref{exm:homeomorphism_not_invariant}, the set $K$ is always a continuum, and hence $u(K)$ is a single point.
\begin{proof}
Let $\rho\in L^2(\Omega)$ be an upper gradient of $u_\Omega$. Define $\widetilde{\rho}\colon \mathbb{D} \to [0, \infty)$ by 
\[
\widetilde{\rho}(z)=\begin{cases}
\rho(z) & \text{ if } \ z \in \Omega\\
0 & \text{ if } \ z \in K\cap \mathbb{D}
\end{cases}.\]
We claim that $\widetilde{\rho}$ is an upper gradient of $u$. By definition, for any absolutely continuous path $\gamma\subset \Omega$ with end points $x$ and $y$ one has
\[
\ell(u\circ \gamma)\leq \int_\gamma \rho=\int_\gamma \widetilde{\rho}.
\]
By the continuity of $u$, the same upper gradient inequality holds if the interior of $\gamma$ lies in $\Omega$ but $x$ and $y$ potentially lie in $K \cap \mathbb{D}$. Now let $\gamma\colon [a,b]\to \mathbb{D}$ be an absolutely continuous path between the points $x,y\in \mathbb{D}$. Then $\gamma^{-1}(\Omega)\subset[a,b]$ is a countable union of relatively open subintervals $I_i$. Thus, since $\mathcal{H}^1(u(K))=0$, one has
\[
d(u(x),u(y))\leq \ell(u\circ \gamma)=\sum_i \ell(u \circ \gamma|_{\bar{I}_i})\leq \sum_i \int_{\gamma|_{\bar{I}_i}}\widetilde{\rho}\leq \int_\gamma \widetilde{\rho}.
\]
We have shown that $\widetilde{\rho}$ is an upper gradient of $u$ and hence $u\in N^{1,2}(\mathbb{D},X)$ and $E^2_+(u)\leq E^2_+(u|_\Omega)$. The inequality $E^2_+(u_{|\Omega})\leq E^2_+(u)$ holds because every upper gradient of $u$ restricts to an upper gradient of $u|_\Omega$.
\end{proof}

\subsection{Fibers of energy-minimizing parametrizations} \label{sec:other_continua}

 In this section, we give the proofs of \Cref{thm:other_continua} and \Cref{thm:uniqueness}.
 \begin{proof}[Proof of \Cref{thm:other_continua}]
 Let $z \in Z$ and $u\colon \overline{\mathbb{D}} \to Z$ be as in \Cref{thm:other_continua}. Let $F \subset \overline{\mathbb{D}}$ be the nondegenerate fiber $u^{-1}(z)$. Then $F$ is cell-like.
 
 First, assume that $z\in \partial Z$. Then $F$ intersects $\mathbb{S}^1$ in a proper connected subset. We may also assume without loss of generality that $E\subset \overline{\mathbb{D}}$ intersects $\mathbb{S}^1$ in a single point. To verify this, let $\overline{B}\subset \mathbb{R}^2$ be a closed ball of least radius containing $E$ and let $p\in \partial B\cap E$. Let $\overline{B}'\subset \overline{\mathbb{D}}$ be a closed ball intersecting $\mathbb{S}^1$ in a single point $p'$. The assumption may be achieved by replacing $E$ with $\varphi(E)$, where $\varphi \colon \overline{B}\to \overline{B}'$ is a homeomorphism mapping~$p$ to~$p'$.
 
 We observe that $\Omega_1=\mathbb{D}\setminus E$ and $\Omega_2=\mathbb{D}\setminus F$ are simply connected domains. By the Riemann mapping theorem, there is a conformal diffeomorphism $f_i\colon\mathbb{D} \to\Omega_i$ for each $i \in \{1,2\}$. By Corollary~2.17 and Proposition~2.14 in \cite{Pom:92}, these extend to homeomorphisms $f_1\colon \mathbb{D}\cup I_1 \to \overline{\mathbb{D}}\setminus E$ and $f_1\colon \mathbb{D}\cup I_2 \to \overline{\mathbb{D}}\setminus F$ where $I_1$ and $I_2 $ are open subintervals of~$\mathbb{S}^1$. Let $h\colon \overline{\mathbb{D}}\to \overline{\mathbb{D}}$ be a Möbius transformation such that $h(I_1)=I_2$.
 Define $v\colon \mathbb{D} \to Z$ by
 \[
 v(w)=\begin{cases}
 z & \text{if } w \in E\\
(u\circ f_2 \circ h \circ f_1^{-1})(w) & \text{if } w \notin E
 \end{cases}.
 \]
 Note that $v$ is continuous, $v|_{\mathbb{S}^1}$ is a monotone parametrization of $\Gamma$ and $v^{-1}(z)=E$. Also, by Lemma~\ref{lemm:surg} and conformal invariance of energy, one has $v\in N^{1,2}(\mathbb{D},Z)$ and
 \[
 E^2_+(v)=E^2_+(v|_{\Omega_1})=E^2_+(u|_{\Omega_2})\leq E^2_+(u).
 \]
 Thus $v$ is an energy-minimizing parametrization.
 
 For the second case, assume that $z\in Z\setminus \partial Z$. Then $\widehat{\mathbb{C}}\setminus F$ and $\widehat{\mathbb{C}}\setminus E$ are simply connected domains containing~$\mathbb{S}^1$. By the Riemann mapping theorem, there is a conformal diffeomorphism $f\colon \widehat{\mathbb{C}}\setminus E \to \widehat{\mathbb{C}}\setminus F$. Let $\Omega \subset \widehat{\C}$ be the complementary domain of the Jordan curve $\Gamma=f^{-1}(\mathbb{S}^1)$ that contains $E$. Let $g\colon \mathbb{D} \to \Omega$ be a conformal diffeomorphism, which again exists by the Riemann mapping theorem. Note that $g$ extends to a homeomorphism $\overline{\mathbb{D}}\to \overline{\Omega}$ by Carathéodory's theorem. Define $v\colon \mathbb{D} \to Z$ by
 \[
 v(w)=\begin{cases}
 z & \text{if } w \in g^{-1}(E)\\
(u\circ f\circ g)(w) & \text{if } w \notin g^{-1}(E)
 \end{cases}.
 \]
 Note that $v$ is continuous, $v|_{\mathbb{S}^1}$ is a monotone parametrization of $\Gamma$ and $v^{-1}(z)=g^{-1}(E)$ is homeomorphic to~$E$. Again, Lemma~\ref{lemm:surg} and the conformal invariance of energy imply that $v\in N^{1,2}(\mathbb{D},Z)$ and
 \[
 E^2_+(v)=E^2_+(v|_{\mathbb{D}\setminus g^{-1}(E)})=E^2_+(u|_{\mathbb{D}\setminus F})\leq E^2_+(u).
 \]
 Thus $v$ is an energy-minimizing parametrization.
 \end{proof}
 \Cref{thm:uniqueness} follows easily from \Cref{thm:other_continua}.
 \begin{proof}[Proof of \Cref{thm:uniqueness}]
 Assume $u$ is an energy-minimizing parametrization of $Z$ and $z\in Z$ is such that $u^{-1}(z)$ is nondegenerate. Then, by Theorem~\ref{thm:other_continua}, for any cell-like subset $E\subset \mathbb{R}^2$ there is an energy-minimizing parametrization $v_E\colon\overline{\mathbb{D}}\to Z$ such that $v_E^{-1}(z)$ is homeomorphic to~$E$. On the other hand, whenever $v_E$ and $v_F$ are equivalent up to conformal diffeomorphism, then~$E$ is homeomorphic to~$F$. Thus it suffices to note that there are infinitely many (in fact, even uncountably many) cell-like subsets of $\mathbb{R}^2$ that are pairwise not homeomorphic.
 \end{proof}
\subsection{A counterexample}

We conclude with an example to show that a situation may happen where one energy minimizer is a homeomorphism while another energy minimizer is not. We observe the following topological fact: if $E, \widehat{E} \subset \widehat{\mathbb{C}}$ are closed and non-separating and $f\colon \widehat{\mathbb{C}} \setminus E \to \widehat{\mathbb{C}} \setminus \widehat{E}$ is a homeomorphism, then $f$ induces a homeomorphism $\widehat{f} \colon \widehat{\mathbb{C}}/\mathcal{E} \to \widehat{\mathbb{C}}/\widehat{\mathcal{E}}$, where $\mathcal{E}$ (resp. $\mathcal{\widehat{E}}$) is the decomposition of $\widehat{\mathbb{C}}$ formed by the components of $E$ (resp. $\widehat{E}$).  More precisely, $\mathcal{E}$ is the partition of $\widehat{\mathbb{C}}$ consisting of the connected components 
of $E$ and the single-point sets $\{x\}$, $x \in \widehat{\mathbb{C}} \setminus E$, and $\widehat{\mathbb{C}}/\mathcal{E}$ is the corresponding quotient space. See the proof of Lemma~5.2 in \cite{IR:20}.

A compact set $E \subset \mathbb{C}$ is \emph{removable for conformal mappings} if every conformal map $f\colon \widehat{\mathbb{C}} \setminus E \to \Omega$, where $\Omega$ is a domain in the Riemann sphere $\widehat{\mathbb{C}}$, extends to a conformal map $F\colon \widehat{\mathbb{C}} \to \widehat{\mathbb{C}}$, that is, a Möbius transformation. See \cite{You:15} for an overview of the topic. 
\begin{exm} \label{exm:homeomorphism_not_invariant}
Let $E \subset [0,1/2] \times \{0\}$ be a totally disconnected set that is not removable for conformal maps. The existence of such a set follows from Theorem 14 in \cite{AB:50}; 
see Example~6.1 and the preceding discussion in \cite{IR:20} for a particular example.  
Let $Z$ be the metric disk obtained by equipping $\overline{\mathbb{D}}$ with the singular Riemannian metric $\chi_{\overline{\mathbb{D}}\setminus E}\cdot g_{||.||}$. Since $E \subset [0,1/2] \times \{0\}$, we see that $Z$ satisfies a $(2\pi)^{-1}$-isoperimetric inequality. Moreover, the natural projection map $P\colon \overline{\mathbb{D}} \to Z$ is a homeomorphic energy minimizer for $Z$. We also make the observation that $P(E)$ has $\mathcal{H}^1$-measure zero in~$Z$.

Since $E$ is totally disconnected, we may identify $\widehat{\mathbb{C}}$ and $\widehat{\mathbb{C}}/\mathcal{E}$. By the assumption that $E$ is not removable for conformal maps, there exists a set $\widehat{E} \subset \widehat{\mathbb{C}}$ and a conformal map $f\colon \widehat{\mathbb{C}} \setminus \widehat{E} \to \widehat{\mathbb{C}} \setminus E$ that does not extend to a Möbius transformation. Then $\widehat{E}$ cannot be totally disconnected. Otherwise, by our topological observation, $f$ would extend to a homeomorphism $\widehat{f}\colon \mathbb{C}\to \mathbb{C}$ which, by Corollary~5.2 in \cite{You:15}, would be conformal. Let $\widehat{\pi}\colon\widehat{\mathbb{C}}\to \widehat{\mathbb{C}}/\mathcal{\widehat{E}}$ be the canonical projection and $\widehat{f}\colon\widehat{\mathbb{C}}/\mathcal{\widehat{E}}\to\widehat{\mathbb{C}}$ be the homeomorphism induced by~$f$. 
Let $\overline{\Omega}\subset \widehat{\mathbb{C}}$ be the Jordan domain such that $(\widehat{f}\circ\widehat{\pi})(\overline{\Omega})=\overline{\mathbb{D}}$ and $h\colon\overline{\mathbb{D}}\to\overline{\Omega}$ be a conformal map.
Then, as in the previous subsection, we may apply \Cref{lemm:surg} to $K=h^{-1}(\widehat{E})$ to see that $\widehat{P}=P \circ \widehat{f}\circ \widehat{\pi}\circ h$ defines an energy-minimizing parametrization. On the other hand, since $\widehat{\pi}$ fails to be injective on $\Omega$, $\widehat{P}$ is not a homeomorphism.
\end{exm}

\subsection*{Acknowledgments}
We thank Alexander Lytchak for bringing to our attention Questions~1.4 and~1.5 in \cite{LW:18a} and for many helpful conversations. 
We also benefited from Malik Younsi's answer to a question we asked on MathOverflow, accessible at https://mathoverflow.net/questions/350593/, which became the basis for \Cref{exm:homeomorphism_not_invariant} below, and further discussions with him about removable sets. Finally, we thank the referee for a careful reading of the paper and useful feedback.

\bibliographystyle{abbrv}
\bibliography{bibliography}

\begin{thebibliography}{HKST15}

\bibitem[AB50]{AB:50}
Lars Ahlfors and Arne Beurling.
\newblock Conformal invariants and function-theoretic null-sets.
\newblock {\em Acta Math.}, 83:101--129, 1950.

\bibitem[BBI01]{BBI:01}
Dmitri Burago, Yuri Burago, and Sergei Ivanov.
\newblock {\em A course in metric geometry}, volume~33 of {\em Graduate Studies
  in Mathematics}.
\newblock American Mathematical Society, Providence, RI, 2001.

\bibitem[BK02]{BK:02}
Mario Bonk and Bruce Kleiner.
\newblock Quasisymmetric parametrizations of two-dimensional metric spheres.
\newblock {\em Invent. Math.}, 150(1):127--183, 2002.

\bibitem[{Cre}20]{Cre:20}
Paul {Creutz}.
\newblock {Rigidity of the Pu inequality and quadratic isoperimetric constants
  of normed spaces}.
\newblock {\em arXiv e-prints}, page arXiv:2004.01076, 2020.

\bibitem[Cre21]{Cre:19}
Paul Creutz.
\newblock Space of minimal discs and its compactification.
\newblock {\em Geom. Dedicata}, 210:151--164, 2021.

\bibitem[Crear]{Cre:ar}
Paul Creutz.
\newblock Plateau's problem for singular curves.
\newblock {\em Comm. Anal. Geom.}, to appear.

\bibitem[Dav86]{Dav:86}
Robert~J. Daverman.
\newblock {\em Decompositions of manifolds}, volume 124 of {\em Pure and
  Applied Mathematics}.
\newblock Academic Press, Inc., Orlando, FL, 1986.

\bibitem[Dou31]{Dou:31}
Jesse Douglas.
\newblock Solution of the problem of {P}lateau.
\newblock {\em Trans. Amer. Math. Soc.}, 33(1):263--321, 1931.

\bibitem[Fed65]{Fed:65}
Herbert Federer.
\newblock Some theorems on integral currents.
\newblock {\em Trans. Amer. Math. Soc.}, 117:43--67, 1965.

\bibitem[Fed69]{Fed:69}
Herbert Federer.
\newblock {\em Geometric measure theory}.
\newblock Die Grundlehren der mathematischen Wissenschaften, Band 153.
  Springer-Verlag New York Inc., New York, 1969.

\bibitem[FM13]{FM:13}
A.~Figalli and F.~Maggi.
\newblock On the isoperimetric problem for radial log-convex densities.
\newblock {\em Calc. Var. Partial Differential Equations}, 48(3-4):447--489,
  2013.

\bibitem[Fre92]{Fre:92}
D.~H. Fremlin.
\newblock Spaces of finite length.
\newblock {\em Proc. London Math. Soc. (3)}, 64(3):449--486, 1992.

\bibitem[GOR73]{GOR:73}
R.~D. Gulliver, II, R.~Osserman, and H.~L. Royden.
\newblock A theory of branched immersions of surfaces.
\newblock {\em Amer. J. Math.}, 95:750--812, 1973.

\bibitem[HKST01]{HKST:01}
Juha Heinonen, Pekka Koskela, Nageswari Shanmugalingam, and Jeremy~T. Tyson.
\newblock Sobolev classes of {B}anach space-valued functions and quasiconformal
  mappings.
\newblock {\em J. Anal. Math.}, 85:87--139, 2001.

\bibitem[HKST15]{HKST:15}
Juha Heinonen, Pekka Koskela, Nageswari Shanmugalingam, and Jeremy~T. Tyson.
\newblock {\em Sobolev spaces on metric measure spaces}, volume~27 of {\em New
  Mathematical Monographs}.
\newblock Cambridge University Press, Cambridge, 2015.
\newblock An approach based on upper gradients.

\bibitem[HW41]{HW:41}
Witold Hurewicz and Henry Wallman.
\newblock {\em Dimension {T}heory}.
\newblock Princeton Mathematical Series, v. 4. Princeton University Press,
  Princeton, N. J., 1941.

\bibitem[{Iko}19]{Iko:19}
Toni {Ikonen}.
\newblock {Uniformization Of Metric Surfaces Using Isothermal Coordinates}.
\newblock {\em arXiv e-prints}, page arXiv:1909.09113, 2019.

\bibitem[IR20]{IR:20}
Toni {Ikonen} and Matthew {Romney}.
\newblock {Quasiconformal geometry and removable sets for conformal mappings}.
\newblock {\em arXiv e-prints}, page arXiv:2006.02776, 2020.

\bibitem[Jos94]{Jos:94}
J\"{u}rgen Jost.
\newblock Equilibrium maps between metric spaces.
\newblock {\em Calc. Var. Partial Differential Equations}, 2(2):173--204, 1994.

\bibitem[Jos17]{Jos:17}
J\"{u}rgen Jost.
\newblock {\em Riemannian geometry and geometric analysis}.
\newblock Universitext. Springer, Cham, seventh edition, 2017.

\bibitem[Kar07]{Kar:07}
M.~B. Karmanova.
\newblock Area and co-area formulas for mappings of the {S}obolev classes with
  values in a metric space.
\newblock {\em Sibirsk. Mat. Zh.}, 48(4):778--788, 2007.

\bibitem[LW17a]{LW:17}
Alexander Lytchak and Stefan Wenger.
\newblock Area minimizing discs in metric spaces.
\newblock {\em Arch. Ration. Mech. Anal.}, 223(3):1123--1182, 2017.

\bibitem[LW17b]{LW:17b}
Alexander Lytchak and Stefan Wenger.
\newblock Energy and area minimizers in metric spaces.
\newblock {\em Adv. Calc. Var.}, 10(4):407--421, 2017.

\bibitem[LW18a]{LW:18a}
Alexander Lytchak and Stefan Wenger.
\newblock Intrinsic structure of minimal discs in metric spaces.
\newblock {\em Geom. Topol.}, 22(1):591--644, 2018.

\bibitem[LW18b]{LW:18b}
Alexander Lytchak and Stefan Wenger.
\newblock Isoperimetric characterization of upper curvature bounds.
\newblock {\em Acta Math.}, 221(1):159--202, 2018.

\bibitem[LW20]{LW:20}
Alexander Lytchak and Stefan Wenger.
\newblock Canonical parameterizations of metric disks.
\newblock {\em Duke Math. J.}, 169(4):761--797, 2020.

\bibitem[Mes00]{M:00}
Chikako Mese.
\newblock Some properties of minimal surfaces in singular spaces.
\newblock {\em Trans. Amer. Math. Soc.}, 352(9):3957--3969, 2000.

\bibitem[MHH11]{MHH:11}
Frank Morgan, Sean Howe, and Nate Harman.
\newblock Steiner and {S}chwarz symmetrization in warped products and fiber
  bundles with density.
\newblock {\em Rev. Mat. Iberoam.}, 27(3):909--918, 2011.

\bibitem[Mor48]{Mor:49}
Charles~B. Morrey, Jr.
\newblock The problem of {P}lateau on a {R}iemannian manifold.
\newblock {\em Ann. of Math. (2)}, 49:807--851, 1948.

\bibitem[MR02]{MR:02}
Frank Morgan and Manuel Ritor\'{e}.
\newblock Isoperimetric regions in cones.
\newblock {\em Trans. Amer. Math. Soc.}, 354(6):2327--2339, 2002.

\bibitem[MZ10]{MZ:10}
Chikako Mese and Patrick~R. Zulkowski.
\newblock The {P}lateau problem in {A}lexandrov spaces.
\newblock {\em J. Differential Geom.}, 85(2):315--356, 2010.

\bibitem[Nik79]{Nik79}
I.~G. Nikolaev.
\newblock Solution of the {P}lateau problem in spaces of curvature at most
  {$K$}.
\newblock {\em Sibirsk. Mat. Zh.}, 20(2):345--353, 459, 1979.

\bibitem[Oss70]{Oss:70}
Robert Osserman.
\newblock A proof of the regularity everywhere of the classical solution to
  {P}lateau's problem.
\newblock {\em Ann. of Math. (2)}, 91:550--569, 1970.

\bibitem[OvdM14]{OvdM:14}
Patrick Overath and Heiko von~der Mosel.
\newblock Plateau's problem in {F}insler 3-space.
\newblock {\em Manuscripta Math.}, 143(3-4):273--316, 2014.

\bibitem[Pom92]{Pom:92}
Ch. Pommerenke.
\newblock {\em Boundary behaviour of conformal maps}, volume 299 of {\em
  Grundlehren der Mathematischen Wissenschaften [Fundamental Principles of
  Mathematical Sciences]}.
\newblock Springer-Verlag, Berlin, 1992.

\bibitem[PvdM17]{PvdM:17}
Sven Pistre and Heiko von~der Mosel.
\newblock The {P}lateau problem for the {B}usemann-{H}ausdorff area in
  arbitrary codimension.
\newblock {\em Eur. J. Math.}, 3(4):953--973, 2017.

\bibitem[Rad30]{Rad:30}
Tibor Rad\'{o}.
\newblock On {P}lateau's problem.
\newblock {\em Ann. of Math. (2)}, 31(3):457--469, 1930.

\bibitem[Raj17]{Raj:17}
Kai Rajala.
\newblock Uniformization of two-dimensional metric surfaces.
\newblock {\em Invent. Math.}, 207(3):1301--1375, 2017.

\bibitem[Res68]{Res:68}
Ju.~G. Reshetnyak.
\newblock Non-expansive maps in a space of curvature no greater than {$K$}.
\newblock {\em Sibirsk. Mat. \v Z.}, 9:918--927, 1968.

\bibitem[Res93]{Res:93}
Yu.~G. Reshetnyak.
\newblock Two-dimensional manifolds of bounded curvature.
\newblock In {\em Geometry, {IV}}, volume~70 of {\em Encyclopaedia Math. Sci.},
  pages 3--163, 245--250. Springer, Berlin, 1993.

\bibitem[Rom19]{Rom:19}
Matthew Romney.
\newblock Singular quasisymmetric mappings in dimensions two and greater.
\newblock {\em Adv. Math.}, 351:479--494, 2019.

\bibitem[RRR19]{RRR:19}
Kai {Rajala}, Martti {Rasimus}, and Matthew {Romney}.
\newblock {Uniformization with infinitesimally metric measures}.
\newblock {\em arXiv e-prints}, page arXiv:1907.07124, 2019.

\bibitem[Spa81]{Spa:81}
Edwin~H. Spanier.
\newblock {\em Algebraic topology}.
\newblock Springer-Verlag, New York-Berlin, 1981.
\newblock Corrected reprint.

\bibitem[Staar]{Sta:ar}
Stephan Stadler.
\newblock {The structure of minimal surfaces in CAT(0) spaces}.
\newblock {\em J. Eur. Math. Soc.}, to appear.
\newblock preprint arXiv:1808.06410.

\bibitem[V{\"{a}}i71]{Vai:71}
Jussi V{\"{a}}is\"{a}l\"{a}.
\newblock {\em Lectures on {$n$}-dimensional quasiconformal mappings}.
\newblock Lecture Notes in Mathematics, Vol. 229. Springer-Verlag, Berlin-New
  York, 1971.

\bibitem[Wil49]{Wil:49}
Raymond~Louis Wilder.
\newblock {\em Topology of {M}anifolds}.
\newblock American Mathematical Society Colloquium Publications, vol. 32.
  American Mathematical Society, New York, N. Y., 1949.

\bibitem[Wil08]{Wil:08}
Kevin Wildrick.
\newblock Quasisymmetric parametrizations of two-dimensional metric planes.
\newblock {\em Proc. Lond. Math. Soc. (3)}, 97(3):783--812, 2008.

\bibitem[Wil10]{Wil:10}
Kevin Wildrick.
\newblock Quasisymmetric structures on surfaces.
\newblock {\em Trans. Amer. Math. Soc.}, 362(2):623--659, 2010.

\bibitem[You15]{You:15}
Malik Younsi.
\newblock On removable sets for holomorphic functions.
\newblock {\em EMS Surv. Math. Sci.}, 2(2):219--254, 2015.

\end{thebibliography}

\end{document}